\documentclass[12pt, preprint]{elsarticle}
\usepackage{amsthm,amssymb,amsfonts,amsmath,lineno,mathtools,graphicx}

\usepackage{epsfig,xspace,algorithm,verbatim}
\usepackage[noend]{algpseudocode}
\usepackage{color}
\usepackage{enumerate}
\usepackage{array}
\usepackage{tikz}
\usepackage {csquotes} 
\usepackage{caption}
\usepackage{grffile}

\newtheorem{theorem}{Theorem}[section]
\newtheorem{lemma}[theorem]{Lemma}
\newtheorem{proposition}[theorem]{Proposition}

\newtheorem{definition}[theorem]{Definition}
\newtheorem{remark}[theorem]{Remark}

\newcommand{\bdf}{\begin{df} \begin{rm}}
\newcommand{\edf}{\end{rm} \end{df}}

\makeatletter
\DeclareRobustCommand{\cev}[1]{%
  \mathpalette\do@cev{#1}%
}
\newcommand{\do@cev}[2]{%
  \fix@cev{#1}{+}%
  \reflectbox{$\m@th#1\vec{\reflectbox{$\fix@cev{#1}{-}\m@th#1#2\fix@cev{#1}{+}$}}$}%
  \fix@cev{#1}{-}%
}
\newcommand{\fix@cev}[2]{%
  \ifx#1\displaystyle
    \mkern#23mu
  \else
    \ifx#1\textstyle
      \mkern#23mu
    \else
      \ifx#1\scriptstyle
        \mkern#22mu
      \else
        \mkern#22mu
      \fi
    \fi
  \fi
}

\makeatother

\begin{document}
\begin{frontmatter}

\title{Gracefulness of two nested cycles: a first approach}

\author{Miguel Licona}
\ead{eliconav23@xanum.uam.mx}

\author{Joaqu\'in Tey\corref{cor1}}

\cortext[cor1]{Corresponding author}
\ead{jtey@xanum.uam.mx}

\address{Universidad Aut\'onoma Metropolitana, Iztapalapa\\
         Av. San Rafael Atlixco 186, col.~Vicentina\\
         09340 Iztapalapa, M\'{e}xico}
         

\begin{abstract}

It is known that if a plane graph admits a graceful (resp. near-graceful) labeling, then its semidual admits a conservative (resp. near-conservative) labeling.

In this work we prove that the semidual of a plane graph of size $M$ consisting of two nested cycles is conservative if $M \equiv 0,3 \pmod 4$, and it is near-conservative otherwise.

We also show that for a given integer $m_1 \geq 3$, there exists $m^* > m_1$ such that for $m_2 \geq m^*$, if $m_1+m_2 \equiv 0,3 \pmod 4$ (resp. $m_1+m_2 \equiv 1,2 \pmod 4$), then there exists a graceful (resp. near-graceful) plane graph consisting of two nested cycles with sizes $m_1$ and $m_2$, respectively.
\end{abstract}

\end{frontmatter}

\section{Introduction}

We deal with finite simple graphs, the vertex and edge sets of a graph $G$ are denoted by $V(G)$ and $E(G)$, respectively. We denote an edge with end-vertices $v$ and $w$ by $vw\in E(G)$. Let $G$ be a graph and $v\in V(G);$ the {\em degree} of $v$, denoted by $d_G(v)$, is the number of edges incident to $v,$  for short $d(v)$ when the underlying graph is understood.  An {\em internal} vertex is a vertex of degree at least two, we denote by $V_{I}(G)$ the set of all internal vertices of $G.$ We refer the reader to \cite{W} for undefined concepts and terminology.

The interval $[a, b]$ denotes the set of consecutive integers $\{a, a + 1,\ldots, b \},$ where $a \leq b.$ 

For our purpose, a {\em labeling} $\phi$ of a graph $G$ is a bijection from $E(G)$ to a set of positive integers of size $\vert E(G) \vert$. For $G' \subseteq G,$ $\phi_{G'}$ denotes the restriction of $\phi$ to $E(G').$ 

A graph $G$ of size $M$ is {\em graceful} if there exists an injection  $f:V(G)\longrightarrow \{0,\ldots,M\}$ such that $\{|f(u)-f(v)|\}_{(u,v)\in E(G)}=\{1,2,\ldots,M\}.$  Similarly, $G$ is near-graceful if it is not graceful and $\{|f(u)-f(v)|\}_{(u,v)\in E(G)}=\{1,2,\ldots,M-1\}\cup  \{M+1\}.$ 

The notion of graceful graph was introduced by Rosa \cite{R} in the late sixties. Since then, gracefulness on graphs has been one of the most studied graph labeling problem. 

Given an orientation $\overrightarrow{G}$ of a graph $G$, we denote by $A(\overrightarrow{G})$ the set of oriented edges of $\overrightarrow{G}$ ({\em arcs} for short) and either  $\overrightarrow{vw}\in A(\overrightarrow{G})$ or  $\overrightarrow{wv}\in A(\overrightarrow{G})$ if $vw \in E(G).$ 
 
Let $\overrightarrow{G}$ be an orientation of $G$ and $\phi$ a labeling of $G.$ The {\em vertex-sum} of a vertex $u \in V(\overrightarrow{G}),$ denoted by $s_{\overrightarrow{G},\phi}(u),$ is defined as
$$ s_{\overrightarrow{G},\phi}(u)=\sum_{\overrightarrow{vu} \in A(\overrightarrow{G})}\phi(uv) - \sum_{\overrightarrow{uv}\in A(\overrightarrow{G})}\phi(uv).$$ 

When no confusion arises, we simply write $s_{\overrightarrow{G}}(u)$ or $s(u).$

Let $G$ be a graph of size $M$.  A {\em conservative labeling} of $G$ is a pair  $( \overrightarrow{G},\phi),$ where  $\overrightarrow{G}$ is an orientation of $G$ and $\phi$ is a labeling from $E(G)$ to  $[1, M]$ so that at each vertex of degree at least three the vertex-sum is zero. Similarly, $( \overrightarrow{G},\phi)$ is a {\em near-conservative labeling} of $G$ if $\phi(E(G))=[1,M-1]\cup \{M+1\}.$ 

We say that $G$ is {\em conservative} if it admits a conservative labeling. When $G$ is non-conservative and it admits a near-conservative labeling, $G$ is a {\em near-conservative} graph.

The {\em semidual} of a plane graph $G$ is the graph obtained by splitting the associated vertex to the outer face in the dual of $G$, into as many vertices as its degree.

A {\em cycle with chords} is a drawing of a Hamiltonian graph, where all the vertices lie on its {\em base cycle}, the boundary of the outer face. In a {\em plane cycle with chords} the chords do not cross. Plane cycles with chords determine a particular subclass of outerplane graphs. A {\em two nested cycles graph} is a plane cycle with chords whose set of chords induces a cycle. 

In this work we begin to study gracefulness of two nested cycles graphs. It is worth emphasizing that these graphs are Eulerian.

Gracefulness in cycles with chords is a well studied problem. Examples of graceful plane cycles with chords are ladders \cite{AG,M} and several shell-type graphs \cite{DKMTT,EE,KP,LT,L,Ma,MLL}. Graceful cycles with crossing chords are prisms \cite{FG,YW} and M\"obius ladders \cite{G}. It should be noted that all known graceful cycles with chords are non-Eulerian graphs. 

On the other hand, a well known result by Rosa \cite{R} asserts that

\begin{theorem}[\cite{R}]\label{Rosa-C}
	Let $G$ be an Eulerian graph. If $\vert E(G) \vert  \equiv 1,2 \pmod 4$, then $G$ is non-graceful. 
\end{theorem}

So, the ``best we can do" regarding gracefulness of an Eulerian graph $G$ with $\vert E(G) \vert  \equiv 1,2 \pmod 4$ is trying to find a near-graceful labeling for it.

Known Eulerian graceful graphs are cycles of size $M \equiv 0,3 \pmod 4$ \cite{R} and several classes of cyclic snakes \cite{AD,Ba0,Ba1,Gn,Mo,S}. Eulerian near-graceful graphs are cycles of size $M \equiv 1,2 \pmod 4$ \cite{Ba0} and several classes of cyclic snakes \cite{AD,Ba0,Mo}. In \cite{RS} was proved that every connected graph can be embedded as an induced subgraph of an Eulerian graceful graph.

In \cite{LT} was showed that the semidual of a plane cycle with chords is a tree without vertices of degree two. Besides, it is known that

\begin{theorem}(\cite{BBS})\label{Gr-Sm}
 If a plane graph is graceful, then its dual is conservative. 
 \end{theorem}
 
 It is no hard to realize that the proof of Theorem \ref{Gr-Sm} given in \cite{BBS} can be adapted to show that if a plane graph admits a near-graceful labeling, then its dual admits a near-conservative labeling. 
 
Thus, a necessary condition for a plane cycle with chords to be graceful (resp. near-graceful) is that its semidual (certainly a tree) is conservative (resp. near-conservative).  

Given a set $S$ of disjoint stars, the tree obtained by identifying, in a single vertex, a leaf from each star in $S$ will be called a {\em snowflake}.

\begin{figure}[H]
\begin{center}
\includegraphics[scale=.4]{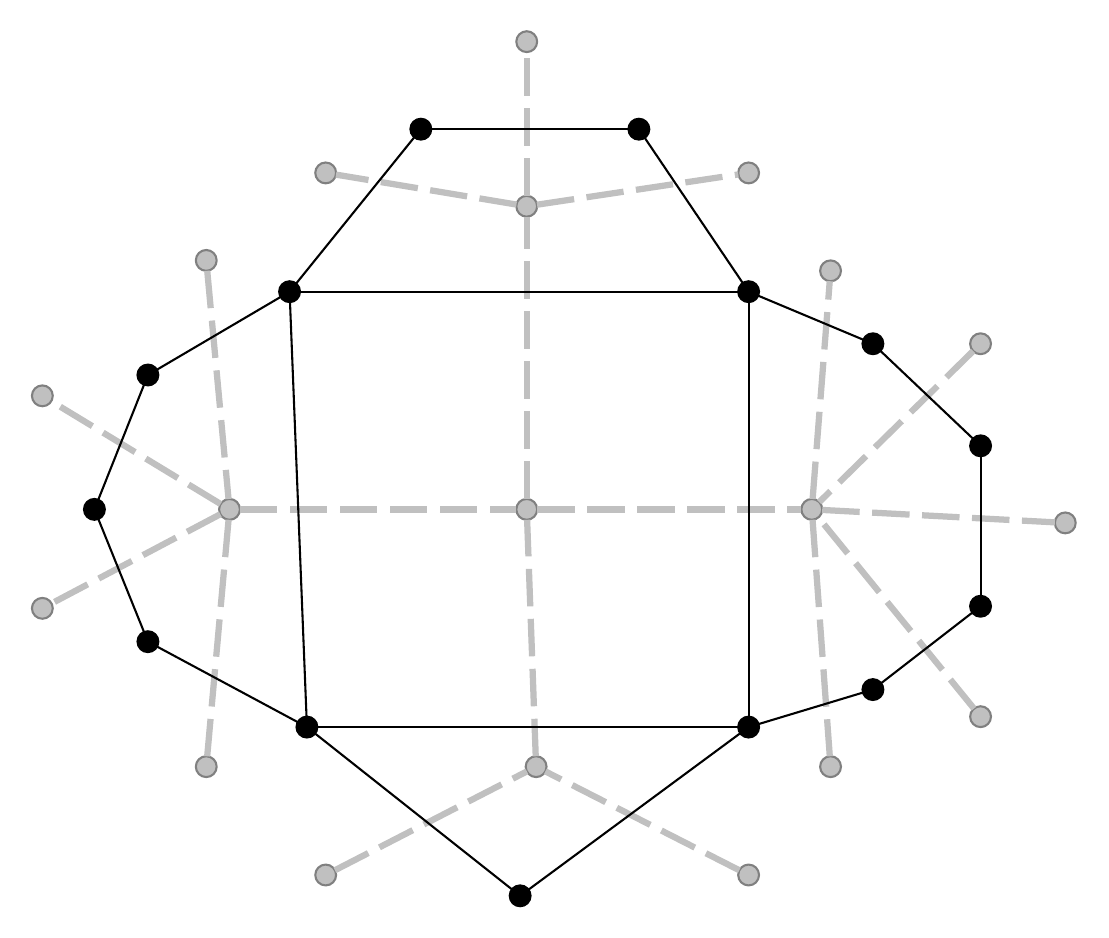} 
\caption{A two nested cycles graph (in black) and its semidual (in gray dashed lines).}
\label{2Nested}
\end{center}
\end{figure}

Notice that the semidual of a two nested cycles graph is a snowflake (see Figure \ref{2Nested}). Hence, gracefulness of two nested cycles graphs can be studied considering two different approaches: the first one by finding graceful (near-graceful) labelings of two nested cycles graphs, and the other one by finding conservative (near-conservative) labelings of snowflakes.

Conservativeness of graphs is an interesting problem on its own. Some conservative graphs are the complete graph \cite{BBS} and several classes of bipartite graphs \cite{ABD,Arch4,BBS,DW}. In \cite{GT} was proved that a galaxy (a disjoint union of stars) of size $M\equiv 0,3 \pmod 4$ is conservative, and it is near-conservative otherwise. Some classes of conservative and near-conservative trees can be found in \cite{LT}.

It should be noted that papers \cite{ABD,Arch4,DW} concern integer Heffter arrays, but as it was pointed out in \cite{GT}, the problem of finding integer Heffter arrays is equivalent to showing that certain biregular bipartite graphs are conservative.

The rest of the text is organized as follows. In Section \ref{Sect-2} we prove that for a given integer $m_1 \geq 3$, there exists $m^* > m_1$ such that for $m_2 \geq m^*$, if $m_1+m_2 \equiv 0,3 \pmod 4$ (resp. $m_1+m_2 \equiv 1,2 \pmod 4$), then there exists a graceful (resp. near-graceful) two nested cycles graph consisting of two cycles having sizes $m_1$ and $m_2$, respectively. In Section \ref{Sect-3} we introduce definitions and some results that will be helpful to developing Section \ref{Sect-4}. Finally, in Section \ref{Sect-4} we show that  a snowflake of size $M$ is conservative if $M \equiv 0,3 \pmod 4$, and it is near-conservative otherwise.

\section{Graceful and near-graceful two nested cycles}\label{Sect-2}

We denote by $\mathcal{N}_g(m,n)$ (resp. $\mathcal{N}_{n\text{-}g}(m,n)$) the set of graceful (resp. near-graceful) two nested cycles graphs with base cycle $C_n$ and cycle chord $C_m$.

\begin{theorem}\label{Grace2nested}
	Let  $m_1 \geq 3$. For $m_2 \geq m_1(2m_1-1)$ we have
	\begin{itemize}
		\item if $m_1+m_2 \equiv 0,3 \pmod 4$, then $\mathcal{N}_g(m_1,m_2) \neq \emptyset$,
		\item if $m_1+m_2 \equiv 1,2 \pmod 4$, then $\mathcal{N}_{n\text{-}g}(m_1,m_2) \neq \emptyset$.
	\end{itemize}
\end{theorem}

\begin{proof}

Let $m=m_1+m_2$. Then $m \equiv t, 3-t \pmod 4$, where $t \in \{0,1 \}$. We proceed to construct a graph $G$ in $\mathcal{N}_g(m_1,m_2)$ (resp.  $\mathcal{N}_{n\text{-}g}(m_1,m_2)$) if $m \equiv 0,3 \pmod 4$ (resp. $m \equiv 1,2 \pmod 4$). So, $G$ is given by two nested cycles, where the interior cycle $C_{m_1}$ has size $m_1$ and the exterior cycle $C_{m_2}$ has size $m_2$. Let  $V(G)=V(C_{m_2})=\{v_{1},v_{2},\ldots,v_{m_2}\}$, $E(C_{m_2})=\{e_{1},e_{2},\ldots,e_{m_2}\},$ where for $1 \leq i \leq m_2-1, e_i=v_iv_{i+1}$ and $e_{m_2}=v_{m_2}v_1$. For $j \geq 3$, let $s_j= 2\sum_{i=1}^{j-1} i$. The cycle chord of $G$ is defined as $(w_1,w_2, \dots,w_{m_1})$, where
\begin{center}
$w_j:=
	\begin{cases*}
		v_1 &
		\mbox{if $j=1$},\\
		v_{5-2t} &
		\mbox{if $j=2$},\\
		v_{s_j-2t+1} &
		\mbox{if $3 \leq j \leq m_1$}.\\
	\end{cases*}$
\end{center}

Let $c_w=\dfrac{(m_1-1)m_1}{2}$ and $f_w > c_w$. Define a bijection 

$f: E(C_{m_2})  \rightarrow [m_1,m-1] \cup \{m+t\} \setminus \{c_w \}$ as follows.

\begin{center}
$f(e_i):=
	\begin{cases*}
	        m+t &
		\mbox{if $i=1$},\\
		m+1-i &
		\mbox{if $i \in [2, m-f_w]$},\\
		m-i &
		\mbox{if $i \in [m-f_w+1, m-c_w-1]$},\\
		m-1-i &
		\mbox{if $i \in [m-c_w, m_2-1]$ and $m_1 \geq 4$},\\
		f_w &
		\mbox{if $i=m_2$}.\\
	\end{cases*}$
\end{center}

%If $m \equiv 1,2 \pmod 4$, then set $f(e_1):=m+1$.
 
Now $\phi: V(G)  \rightarrow [0,m-1 ] \cup \{m+t \}$ is given by,

\begin{itemize}
	\item $\phi(v_{1}):=0,$
	
	\item 
	for  $2\leq i\leq m_2,$ \hspace{0.2cm} $\phi(v_{i}):=\phi(v_{i-1})+(-1)^{i}f(e_{i-1})$.
	
In the previous definition of $\phi$, the value of $f_w$ varies according to the following cases.

{\em Case 1.} 

$(m \equiv 0 \pmod 4$  and  $m_1 \equiv 0 \pmod 2)$ or $(m \equiv 3 \pmod 4$ and  $m_1 \equiv 1 \pmod 2)$.

\begin{center}
 $f_w:= \dfrac{m+m_1}{2}$.
\end{center}

{\em Case 2.} 

$(m \equiv 2 \pmod 4$ and  $m_1 \equiv 0 \pmod 2)$ or $ (m \equiv 1 \pmod 4$ and  $m_1 \equiv 1 \pmod 2)$.

\begin{center}
 $f_w:= \dfrac{m+m_1}{2}+1$.
\end{center}

{\em Case 3.} 

$(m \equiv 0 \pmod 4$ and  $m_1 \equiv 1 \pmod 2)$ or $ (m \equiv 3 \pmod 4$ and  $m_1 \equiv 0 \pmod 2)$.

\begin{center}
 $f_w:= \dfrac{m-m_1+1}{2}$.
\end{center}

{\em Case 4.} 

$ (m \equiv 2 \pmod 4$ and  $m_1 \equiv 1\pmod 2)$ or $(m \equiv 1 \pmod 4$ and  $m_1 \equiv 0\pmod 2)$.

\begin{center}
 $f_w:= \dfrac{m-m_1+1}{2}+1$.
\end{center}
\end{itemize}

Now we will show that $\phi$ is a graceful labeling of $G$, if Case 1 holds.

In this case $t=0$. Since the sequence $(f(e_1),f(e_2),\dots,f(e_{m_2-1}))$ is strictly decreasing, the sequence $(\phi(v_1),\dots, \phi(v_{m_2}))$ has no repeated elements, so $\phi$ is an injection.
Also $\phi(v_{m_2})= m-(\dfrac{m-m_1}{2})=f_w$ and, by construction, we have $\{ \vert \phi(u)- \phi(v) \vert \}_{uv \in E(C_{m_2})}=\{ f(e_i): 1\leq i \leq m_2 \}=[m_1,m] \setminus \{c_w \}$.

On the other hand,  $\phi(v_{2i-1})=i-1$, for $1 \leq i \leq c_w$. So

\begin{center}
$\phi(w_j):=
	\begin{cases*}
	 0 &
		\mbox{if $j=1$},\\
		j &
		\mbox{if $j \in \{2,3\}$},\\
		\phi(w_{j-1})+j &
		\mbox{if $4 \leq j \leq m_1$}.\\
	\end{cases*}$
\end{center}

Hence, $\{ \vert \phi(u)- \phi(v) \vert \}_{uv \in E(C_{m_1})}= [1,m_1-1] \cup \{c_w\}$. Therefore,
$$\{ \vert \phi(u)- \phi(v) \vert \}_{uv \in E(G)}=\{ \vert \phi(u)- \phi(v) \vert \}_{uv \in E(C_{m_2})} \cup \{ \vert \phi(u)- \phi(v) \vert \}_{uv \in E(C_{m_1})}=[1, m]$$ as needed. 

The analysis of the remaining cases is similar to that of Case 1 and will be omitted. For the sake of clarity we give examples of each cases below.

\end{proof}
The following figures are examples of Theorem \ref{Grace2nested} when $m_1 \equiv 1 \pmod 2.$ Blue squares represent the (labeled) vertices of $V(G)$. In strong blue the vertices of $C_{m_1}$. Green squares contain the labels of the edges of $C_{m_2}$. 

\begin{figure}[H]
\begin{center}
\includegraphics[scale=.55]{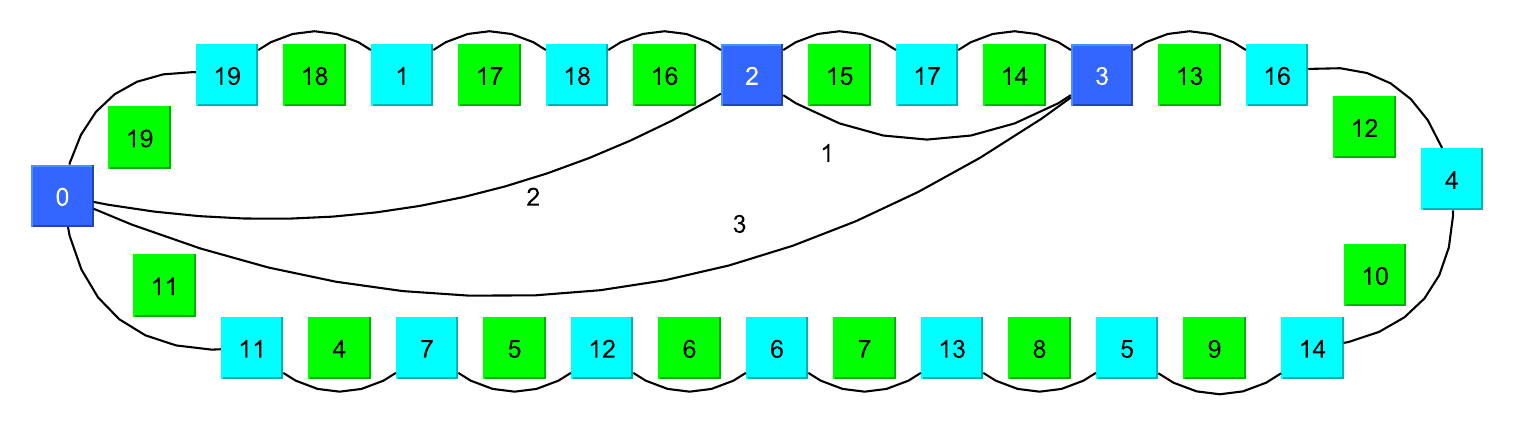} \caption{Example of Case 1. $\mathcal{N}_g(3,16) \neq \emptyset$.}
\label{F1}
\end{center}
\end{figure}

\begin{figure}[H]
\begin{center}
\includegraphics[scale=.53]{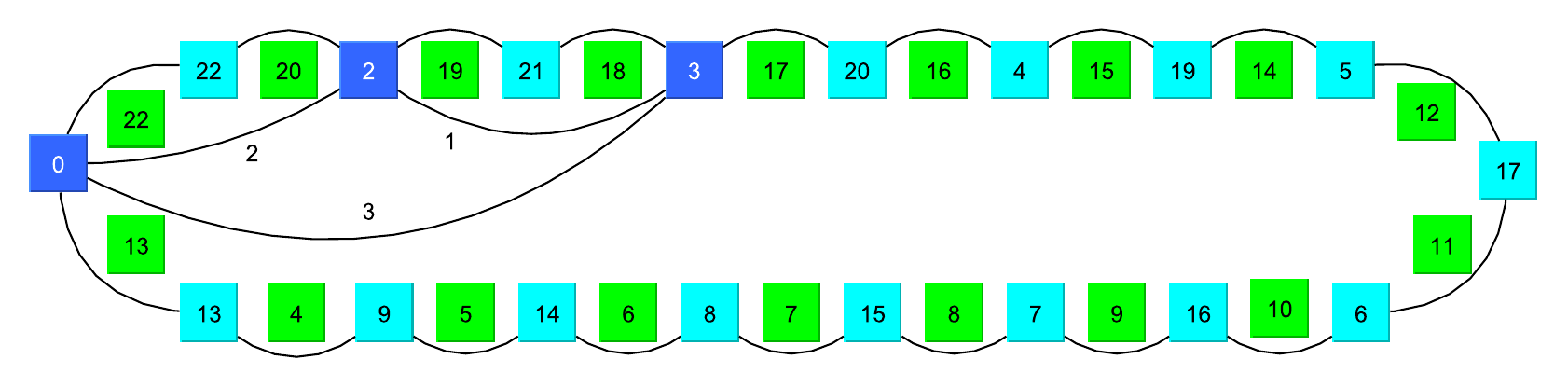} \caption{Example of Case 2. $\mathcal{N}_{n\text{-}g}(3,18) \neq \emptyset$.}
\label{F2}
\end{center}
\end{figure}

\begin{figure}[H]
\begin{center}
\includegraphics[scale=.54]{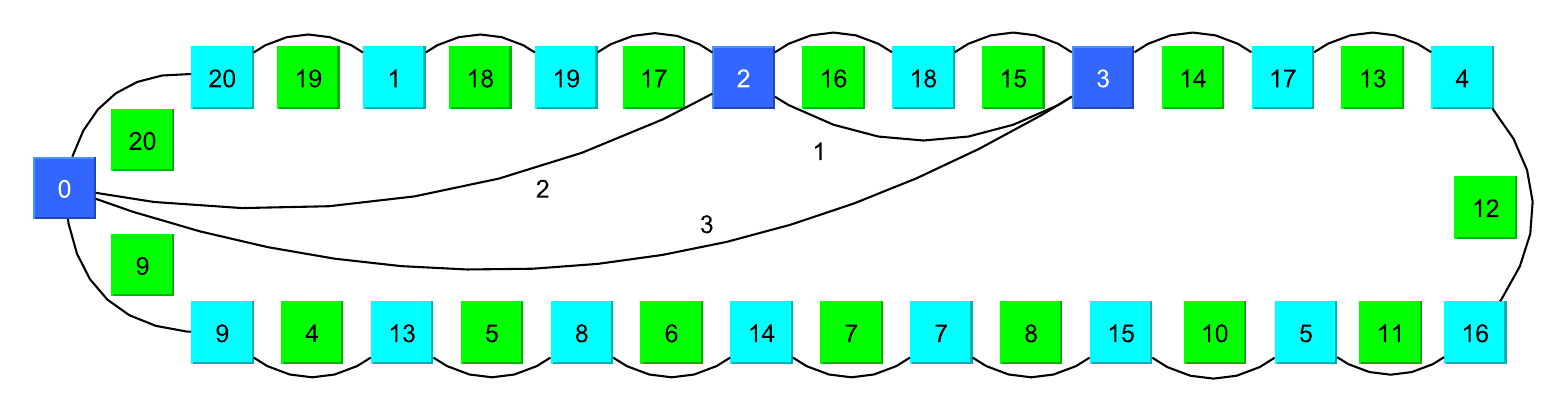} \caption{Example of Case 3. $\mathcal{N}_g(3,17)\neq \emptyset$.}
\label{F3}
\end{center}
\end{figure}

\begin{figure}[H]
\begin{center}
\includegraphics[scale=.54]{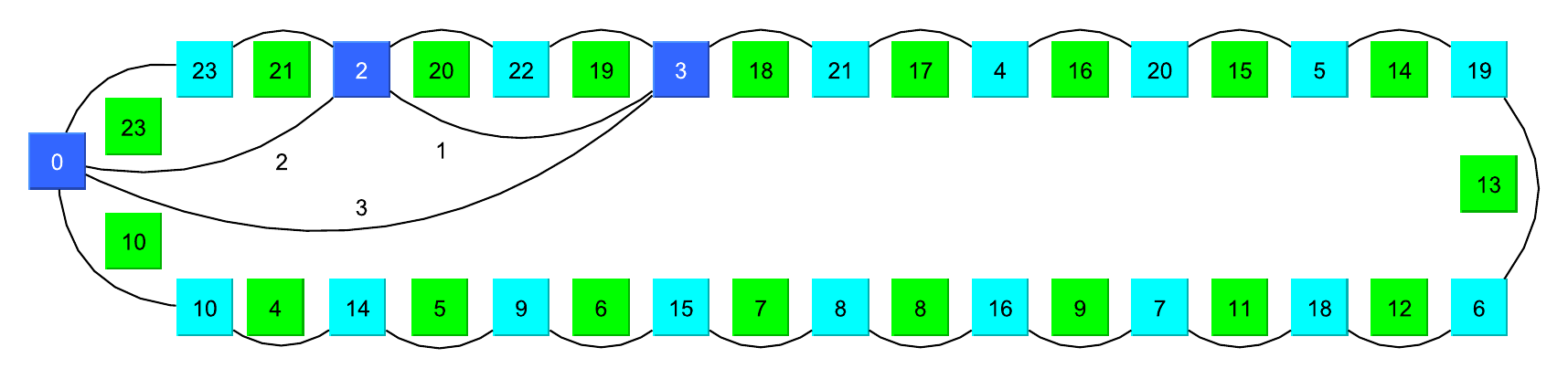} 
\caption{Example of Case 4. $\mathcal{N}_{n\text{-}g}(3,19) \neq \emptyset$.}
\label{F4}
\end{center}
\end{figure}

Now, we will show examples when $m_{1}\equiv 0 \ (mod \ 2),$ more specific when  $m_{1}=4$, and  so $c_w=6$. A number with an upper-bar in the given sequence $(\phi(v_1),\ldots,\phi(v_{m_2}))$ represents a (labeled) vertex of $C_{m_1}$. 

\textit{Example of Case 1.} Let $m_{2}=28$. Then $m=32$, $f_{w}=18$ and  
\begin{center}
$f(e_i):=
	\begin{cases*}
	        32 &
		\mbox{if $i=1$},\\
		33-i &
		\mbox{if $i \in [2, 14]$, \hspace{0.3cm} $(31,30,\ldots,19)$},\\
		32-i &
		\mbox{if $i \in [15, 25]$, \hspace{0.3cm} $(17,16,\ldots,7)$},\\
		31-i &
		\mbox{if $i \in [26, 27]$, \hspace{0.3cm} $(5,4)$},\\
		18 &
		\mbox{if $i=28$}.\\
	\end{cases*}$
\end{center}

$\mathcal{N}_g(4,28)\neq \emptyset$ and $(\phi(v_1),\ldots,\phi(v_{28}))$ is given by
$$(\bar{0},32,1,31,\bar{2},30,\bar{3},29,4,28,5,27,\bar{6},26,7,24,8,23,9,22,10,21,11,20,12,19,14,18).$$ 

\textit{Example of Case 2.} Let $m_{2}=30$. Then $m=34$, $f_{w}=20$ and  
\begin{center}
$f(e_i):=
	\begin{cases*}
	        35 &
		\mbox{if $i=1$},\\
		35-i &
		\mbox{if $i \in [2, 14]$, \hspace{0.3cm} $(33,32,\ldots,21)$},\\
		34-i &
		\mbox{if $i \in [15, 27]$, \hspace{0.3cm} $(19,18,\ldots,7)$},\\
		33-i &
		\mbox{if $i \in [28, 29]$, \hspace{0.3cm} $(5,4)$},\\
		20 &
		\mbox{if $i=30$}.\\
	\end{cases*}$
\end{center}

$\mathcal{N}_{n\text{-}g}(4,30) \neq \emptyset$ and $(\phi(v_1),\ldots,\phi(v_{30}))$ is given by
$$(\bar{0},35,\bar{2},34,\bar{3},33,4,32,5,31,\bar{6},30,7,29,8,27,9,26,10,25,11,24,12,23,13,22,14,21,16,20).$$

\textit{Example of Case 3.} Let $m_{2}=31$. Then $m=35$, $f_{w}=16$ and  
\begin{center}
$f(e_i):=
	\begin{cases*}
	        35 &
		\mbox{if $i=1$},\\
		36-i &
		\mbox{if $i \in [2, 19]$, \hspace{0.3cm} $(34,33,\ldots,17)$},\\
		35-i &
		\mbox{if $i \in [20, 28]$, \hspace{0.3cm} $(15,14,\ldots,7)$},\\
		34-i &
		\mbox{if $i \in [29, 30]$, \hspace{0.3cm} $(5,4)$},\\
		16 &
		\mbox{if $i=31$}.\\
	\end{cases*}$
\end{center}

$\mathcal{N}_g(4,31)\neq \emptyset$ and $(\phi(v_1),\ldots,\phi(v_{31}))$ is given by
$$(\bar{0},35,1,34,\bar{2},33,\bar{3},32,4,31,5,30,\bar{6},29,7,28,8,27,9,26,11,25,12,24,13,23,14,22,15,20,16).$$

\textit{Example of Case 4.} Let $m_{2}=29$. Then $m=33$, $f_{w}=16$ and  
\begin{center}
$f(e_i):=
	\begin{cases*}
	        34 &
		\mbox{if $i=1$},\\
		34-i &
		\mbox{if $i \in [2, 17]$, \hspace{0.3cm} $(32,31,\ldots,17)$},\\
		33-i &
		\mbox{if $i \in [18, 26]$, \hspace{0.3cm} $(15,14,\ldots,7)$},\\
		32-i &
		\mbox{if $i \in [27, 28]$, \hspace{0.3cm} $(5,4)$},\\
		16 &
		\mbox{if $i=29$}.\\
	\end{cases*}$
\end{center}

$\mathcal{N}_{n\text{-}g}(4,29) \neq \emptyset$ and $(\phi(v_1),\ldots,\phi(v_{29}))$ is given by
$$(\bar{0},34,\bar{2},33,\bar{3},32,4,31,5,30,\bar{6},29,7,28,8,27,9,26,11,25,12,24,13,23,14,22,15,20,16).$$ 

\section{Preliminaries for Section 4}\label{Sect-3}

\begin{definition}\label{defik-conser}
Let $G$ be a graph of size $M,$ $t\in \{0,1\},$ $k\geq 0.$  A $k$-$t$-conservative labeling of $G$ is a pair $(\overrightarrow{G},\phi),$ where $\overrightarrow{G}$ is an orientation  of $G$ and $\phi$ is labeling from $E(G)$ to $[1,r]\cup[r+1+k,M+k-1]\cup\{ M+k+t\},$ for $r$ a positive integer, such that at each vertex of degree at least three the vertex-sum is zero.
\end{definition}

Notice that, when $k=0=t,$ $(\overrightarrow{G},\phi)$  is a conservative labeling of $G.$ When $k=0$ and $t=1$ $(\overrightarrow{G},\phi)$ is a near-conservative labeling of $G.$

For $t\in\{0,1\},$  a \textit{$t$-conservative} labeling will refer to a conservative labeling when $t=0$ and a near-conservative labeling when $t=1.$ 

An oriented graph  $\overrightarrow{G}$ is {\em Eulerian} if at each internal vertex of  $\overrightarrow{G}$ the number of the incoming arcs is equal to the number of the outgoing arcs. $( \overrightarrow{G}, \phi)$ is an {\em Eulerian $t$-conservative}  labeling of a graph $G$ if it is a $t$-conservative  labeling of $G$ such that  $\overrightarrow{G}$ is Eulerian. We say that $G$ is {\em Eulerian $t$-conservative} if it admits an Eulerian $t$-conservative labeling. In this section, every labeling of a graph $G$ is associated with an orientation of $G$. So, as abuse of notation, we say that  $\phi$ is an {\em Eulerian labeling} of $G$ if the associated orientation $\overrightarrow{G}$ is Eulerian.

\begin{remark}[\cite{BBS}] \label{fuerteconser}
Let $t\in\{0,1\}$ and $G$ be a graph of size $M.$ If $G$ admits an Eulerian $t$-conservative labeling $( \overrightarrow{G},\phi).$ Then for every integer $\alpha \geq 0,$ the labeling $\phi_{\alpha}$ defined as $\phi_\alpha(e):=\phi(e)+\alpha$ satisfies that $s_{\overrightarrow{G},{\phi_{\alpha}}}(u)=0$ for every vertex $u$ of degree at least three in $G$.
\end{remark}

\begin{remark}[\cite{GT}] \label{star-cons} \label{star}
A star of size $M \equiv 0,3 \ (mod \ 4)$ is conservative, and it is near-conservative otherwise.
\end{remark}

\begin{remark}[\cite{GT}]\label{StarEulerian}
A star of size $M \equiv 0 \ (mod \ 4)$ is Eulerian conservative.
\end{remark}

\begin{definition}[\cite{BBS}]
	Let $T_1$ and $T_2$ be two vertex-disjoint trees and $v_1$ (resp. $v_2$) be a vertex of $T_1$ (resp. $T_2$). $T_1(v_1) \bullet T_2(v_2),$ i.e. the {\em attachment} of $T_1$ and $T_2$ at $(v_1,v_2),$ is the resulting tree of identifying the vertices $v_1$ and $v_2.$
\end{definition}

\begin{proposition}[\cite{LT}]\label{atachment}
	Let $T_1$ and $T_2$ be two vertex-disjoint trees and $v_1$ (resp. $v_2$) an internal vertex of $T_1$ (resp. $T_2$). Then the following statements hold.
	\begin{enumerate}
		\item If $T_1$ is conservative and $T_2$ is Eulerian conservative, then $T_1(v_1) \bullet T_2(v_2)$ is conservative.
		\item If $T_1$ is conservative and $T_2$ admits an Eulerian near-conservative labeling, then $T_1(v_1) \bullet T_2(v_2)$ admits a near-conservative labeling. 
		\item If $T_1$ and $T_2$ are Eulerian conservative, then $T_1(v_1) \bullet T_2(v_2)$ is Eulerian conservative.
		\item If $T_1$ is Eulerian conservative and $T_2$ admits an Eulerian near-conservative labeling, then $T_1(v_1) \bullet T_2(v_2)$ admits an Eulerian near-conservative labeling.
	\end{enumerate}
\end{proposition}
 
For $t \in \{0$, $1\}$, a $t$-\emph{Skolem sequence} of order $n$ is a partition of $[1$, $2n-1]\cup\{2n+t\}$ into a collection of disjoint ordered pairs $\{(a_{i}$, $b_{i})\colon i \in [1$, $n]\}$ such that $a_{i}<b_{i}$ and $b_{i}-a_{i}=i$.

The existence of $t$-Skolem sequences is the key to construct $t$-conservative labelings for odd snowflakes (see Subsection \ref{Sect-4.2}).
 
\begin{theorem}[\cite{ST}, \cite{OK}]
\label{Skolem}
Let $t \in \{0$, $1\}$. There exists a $t$-Skolem sequence of order $n$ if and only if $n\equiv t+1$, $-t\pmod4$.
\end{theorem}

 \textit{Skolem type sequence} will refer to a Skolem sequence or a hooked Skolem sequence which is by the Skolem or O'Keefe method, respectively. Given two sequences $S=(a_{1}, a_{2},\ldots,a_{p})$ and $T=(b_{1},b_{2},\ldots, b_{q})$, $S\ast T$ denotes the sequence $(a_{1}, a_{2},\ldots, a_{p}$, $b_{1}$, $b_{2},\ldots, b_{q})$.
 
\begin{definition}\label{hefftersys}
Let $m_{1},m_{2},\ldots,m_{n}$ be positive integers with $m_{i}\geq 3$ for $1\leq i\leq n,$ $M=\sum\limits_{i=1}^{n} m_{i},$ $k$ be a no negative integer and $t\in \{0,1\}.$
\\
Given a family of $n$ disjoints sequences of integers $S=\{D_{1},\ldots,D_{n}\},$ we say that $S$ is a $(m_{1},m_{2},\ldots,m_{n};k;t)$-Skolem system of order $n$ and size $M$ if it satisfies the following conditions

\begin{enumerate}
\item $|D_{i}|=m_{i}$ for $1\leq i\leq n,$
\item $\sum\limits_{d\in D_{j}} d=0$ for $1\leq i\leq n,$ and
\item $abs(S)$ is a partition of $[1,r]\cup[r+1+k,M+k-1]\cup\{ M+k+t\},$ for a positive integer $r$ (denoted by $r(S)$ for $k>0$), where  $abs(S)$ denotes $\{abs(D_{1}),\ldots,abs(D_{n})\}$ and for $D\in S,$  $abs(D)$  is the set of the absolute values of the elements in $D.$
\end{enumerate} 
\end{definition}  

When no confusion arises, we write $(M; k; t)$-Skolem system of order $n$ instead of $(m_1, m_2,\ldots, m_n; k; t)$-Skolem system of order $n$ and size $M = \sum \limits ^{n}_{i=1} m_i$. Also, if $m_1 = m_2 = \ldots = m_n = m,$ then we write $(m\cdot n; k; t)$-Skolem system.

In what follows, the $(m_1, m_2,\ldots, m_n; k; t)$-Skolem system constructed in \cite{GT} will be refered as the $(m_1, m_2,\ldots, m_n; k; t)$-Skolem system obtained by using the Skolem-O'Keefe Method.

For a given sequence of integers $D=(d_{i}),$ we define $-D=(-d_{i}).$ 

\begin{remark}\label{HS-HS}
If $S$ is an Skolem system and $D \in S$, then $(S \setminus \{D\})\bigcup \{-D\} $  is also an  Skolem system.
\end{remark}

\begin{definition}\label{zerohs}
Let $S=\{D_{1},\ldots,D_{n}\}$  be a $(m_{1},m_{2},\ldots,m_{n};k;t)$-Skolem system. We say that $S$ is a zero-sum  $(m_{1},m_{2},\ldots,m_{n};k;t)$-Skolem system if for $1\leq i\leq n,$ there exists $d_{i}\in D_{i}$ such that $\sum\limits _{i=1}^{n}d_{i}=0.$ We say that $d_{i}$ is the distinguish element of $D_{i}$ and we will denote by $Z_{0}(S)$ the set of distinguish elements.  
\end{definition}

\begin{definition}\label{R-orderMi}
Let $S$ be a $(3\cdot n;k;t)$-Skolem system. We say that an order $(D_{1},\ldots,D_{n})$ of $S$ is an {\em $R$-order} if the following conditions hold.
\begin{itemize}
\item[i.] For every $1\leq i\leq n-1,$ there exist $d\in abs(D_{i})$ and $d'\in abs(D_{i+1})$ such that $|d-d'|=1$, 
\item[ii.] $3n+k+t\in abs(D_{n}).$
\end{itemize}
\end{definition}

\begin{lemma}(\cite{GT})\label{ordenR1}
Let $k,n$ and $t$ be integers such that $k\geq 0,$ $n\geq 1$ and  $t\in \{0,1\}.$ If $n\equiv t+1,-t$ $(mod \ 4),$ then there exists a $(3\cdot n;k;t)$-Skolem system which admits an $R$-order.
\end{lemma}

\begin{remark}\label{ordenRcong3}
Let $k,n,s$ and $t$ be integers such that $k\geq 0,s \geq 2, n=4s+1$ and  $t\in \{0,1\}.$ Then, the $R$-order $(D_{1},\ldots,D_{n})$ of the $(3\cdot n;k;t)$-Skolem, constructed in the proof of Lemma \ref{ordenR1}, holds:
\begin{itemize}
\item[i.] $D_{2s+3}=\{-4, -(n+6s+k),n+6s+k+4\}.$ 
\item[ii.] For every $1\leq i\leq 2s+2$ and $2s+4\leq i\leq 4s+1$ there exist $d\in abs(D_{i})$ and $d'\in abs(D_{i+1})$ such that $|d-d'|=1$,
\item[iii.] For every $2s+3\leq i\leq 4s+1$ there exist $d_{1},d_{2}\in abs(D_{i})$ and $d_{1}',d'_{2}\in abs(D_{i+1})$ such that $d_{1}=d_{1}'+1$ and $d_{2}+1=d_{2}',$ where $d_{1}<d_{2}$ and $d_{1}'<d'_{2}.$
\item[iv.] $3n+k+t\in abs(D_{n}).$
\end{itemize}
\end{remark}

\begin{definition}\label{R'-orderM}
Let $S$ be a $(3 \cdot n;k;t)$-Skolem system. We say that an $R$-order $(D_{1},\ldots,D_{n})$ of $S$ is an {\em $R'$-order} if the following condition holds.
\begin{itemize}
%\itemFor every $1\leq i\leq n-1,$ there exist $d\in abs(D_{i})$ and $d'\in abs(D_{i+1})$ such that $|d-d'|=1$, 
\item[iii.] For some $q$ odd, there exist $d\in abs(D_{q})$ and $d'\in abs(D_{q+1})$  such that  $|d-d'|=2$. 
%\item[iii.] $3n+k+t\in abs(D_{n})$.
\end{itemize}
\end{definition}

\begin{lemma}\label{ordenR2}
Let $k,n$ and $t$ be integers such that $k\geq 0,$ $n\geq 1$ and  $t\in \{0,1\}.$ If $n\equiv 2$ $(mod \ 4),$ then there exists a $(3\cdot n;k;1)$-Skolem system which admits an $R'$-order.

\end{lemma}

\begin{proof}
Let $n=4s+2$ and $S_{n}$ be a Skolem type sequence obtained  using  the following table.
\begin{table}[h!]
\begin{center}
\scalebox{.9}{
\begin{tabular}{| c | c | c | c | c |}
\hline
   & i         &  j       & for $s\geq 2$ and &  we obtain: \\  \hline
1) & $r$       & $4s-r+2$ & $1\leq r\leq 2s$  & $\{(1,4s+1),(2,4s),\ldots,(2s,2s+2)\}$\\
2) & $4s+r+3$  & $8s-r+4$ & $1\leq r\leq s-1$ & $\{(4s+4,8s+3),(4s+5,8s+2),\ldots,(5s+2,7s+5)\}$\\
3) & $5s+r+2$  & $7s-r+3$ & $1\leq r\leq s-1$ & $\{(5s+3,7s+2),(5s+4,7s+1),\ldots,(6s+1,6s+4)\}$\\
4) & $2s+1$    & $6s+2$   &                   & $\{(2s+1,6s+2)\}$\\
5) & $4s+2$    & $6s+3$   &                   & $\{(4s+2,6s+3)\}$ \\ 
6) & $4s+3$    & $8s+5$   &                   & $\{(4s+3,8s+5)\}$\\
7) & $7s+3$    & $7s+4$   &                   & $\{(7s+3,7s+4)\}$\\ \hline
\end{tabular}}
\end{center}
\caption{$1$-Skolem sequence for $n \equiv 2 \ (mod\ 4).$} \label{tabla3}
\end{table}

\begin{remark}\label{ordenparejas}
Let $\{(c_{i},d_{i})\}_{i=1}^{r}$ be a set of pairs obtained in a chosen row of the Table \ref{tabla3}. Observe that, $(c_{1},c_{2},\ldots,c_{r})$ is an increase sequence and $(d_{1},d_{2},\ldots,d_{r})$ is a decrease sequence. 
\end{remark}

In order to obtain a $(3\cdot n;k;1)$-Skolem system $S$ from $S_{n}$, it is require to add  $n+k$  to each number into the pair in $S_n$, i.e, $abs(S)= \bigcup_{i=1}^{n}\{|(a_{i}+n)-(b_{i}+n)|, a_{i}+n+k,b_{i}+n+k\}.$ 

For $s=1$, consider $$\overrightarrow{S}=(\{-4,-(7+k),11+k\},\{-2,-(8+k),10+k\},\{-5,-(9+k),14+k\},$$ $$\{-1,-(16+k),17+k\},\{(-3,-(12+k),15+k\},\{-6,-(13+k),(19+k)\}).$$
\\ 

To make the reminder of this proof easier we will give an order $(D_1,\ldots,D_n)$  to $S_{n}$ and it will induce an $R'$-order in $S$. For $s\geq 2,$  by Remark \ref{ordenparejas} we have the following:
\\

\textit{Case 1.} If $s$ is even, set
\\

$D_{1}:=(6s+3,4s+2),$
\\
$D_{2}:=(1,4s+1),D_{3}:=(2,4s),\ldots,D_{2s}:=(2s,2s+2),$
\\
$D_{2s+1}:=(2s+1,6s+2),$
\\
$D_{2s+2}:=(6s+1,6s+4),D_{2s+3}:=(6s,6s+5)\ldots, D_{3s+1}:=(5s+3,7s+2),$
\\
$D_{3s+2}:=(7s+3,7s+4),$
\\
$D_{3s+3}:=(7s+5,5s+2),\ldots,D_{4s+1}:=(8s+3,4s+4),$
\\
$D_{4s+2}:=(4s+3,8s+5).$
\\

\textit{Case 2.} If $s$ is odd, set 
\\
$D_{1}:=(6s+3,4s+2),$
\\
$D_{2}:=(1,4s+1),D_{3}=(2,4s),\ldots,D_{2s}:=(2s,2s+2),$
\\
$D_{2s+1}:=(2s+1,6s+2),$
\\
$D_{2s+2}:=(6s+1,6s+4),D_{2s+3}:=(6s,6s+5)\ldots,D_{3s}:=(5s+4,7s+1),D_{3s+1}:=(7s+3,7s+4),D_{3s+2}:=(7s+2,5s+3),$
\\
$D_{3s+3}:=(7s+5,5s+2),\ldots,D_{4s+1}:=(8s+3,4s+4),$
\\
$D_{4s+2}:=(4s+3,8s+5).$
\\
So, let $S_{4s+2}:=(D_{1}\ldots,D_{4s+2}).$
\\

In both cases condition \textit{i.} and \textit{ii.} of Definition \ref{R-orderMi} hold. When $s$ is even, consider the sets $\{5s+3+n,7s+2+n\}\in D_{3s+1},$ $\{7s+4+n,7s+3+n\}\in D_{3s+2}$ and $\{7s+5+n,5s+2+n\}\in D_{3s+3}.$ When $s$ is odd, consider the sets  $\{5s+4+n,7s+1+n\}\in D_{3s},$ $\{7s+3+n,7s+4+n\}\in D_{3s+1}$ and $\{7s+2+n,5s+3+n\}\in D_{3s+2}.$  With the previous sets we have condition \textit{iii.} of Definition \ref{R'-orderM}.
\\
Thus,  we  obtain the desired $R'$-order of $S$.
\end{proof}

\begin{lemma}(\cite{GT})\label{3,5hefftersystem}
Let $k, n, n_{1},n_{2}$ and $t$ be integers such that $k \geq 0,$ $n=n_{1}+n_{2},$ $n>0$, and $t\in \{0,1\}$. If $n\equiv t+1,-t \ (mod \ 4)$, then there exists a $(3\cdot n_{1},5\cdot n_{2}; k: t)$-Skolem system $S'$ of order $n$ and size $3\cdot n_{1}+5\cdot n_{2}.$ Moreover, if $n_{2}\equiv 0 \ (mod \ 2)$ and $n_{2}<n_{1}$, there exists $D \in S'$ such that $|D|=3$, and $3\cdot n_{1}+5\cdot n_{2}+k+t\in D_{n}.$
\end{lemma}

\section{Conservative and near-conservative snowflakes}\label{Sect-4}

A {\em galaxy} is a disjoint union of stars.

\begin{definition}\label{defcopo}
Let $G=\bigcup_{i=1}^p S_i$ be a galaxy, where $S_i$ is a star of size $n_i$. The \textit{snowflake} given by $G$, denoted by $\mathbf{C}_{n_{1},\ldots,n_{p}}$, is the tree obtained by identifying,  in a single vertex  $z,$ a leaf from each star in $G$. We say that $z$ is the {\em center} of the snowflake. 
\end{definition}

We say that the snowflake $C_{n_{1},\ldots,n_{p}}$ is \textit{minimum} if  $3\leq n_{i}\leq 6,$ for $1 \leq i \leq p$. We use the notation $C_{3\cdot n_{1},5\cdot n_{2},6\cdot n_{3},4\cdot n_{4}}$ for a minimum snowflake which result for the attachment of $n_1$ stars of size three, $n_2$ stars of size five, $n_3$ stars of size six, and $n_4$ stars of size four.

The following remark is obvious but useful.

\begin{remark}\label{sumapar}
Suppose that $k\geq 2$ and for $1\leq i \leq 2k,$ $n_{i}=1.$ Then, in the sums $n_{1}+\ldots +n_{2k}$, $2+n_{1}+\ldots +n_{2k}$ and $4+n_{1}+\ldots +n_{2k};$ we can change the sign of some $n_{i}'s$ such that the new sum is zero.
\end{remark}

A {\em star-decomposition} of a graph $G$ is a set of stars $\{S_{v_1},S_{v_2},\ldots,S_{v_k}\}$ contained in $G$ such that $\left(E(S_{v_1}),E(S_{v_2}),\ldots,E(S_{v_k})\right)$ is a partition of $E(G)$ and $d_G(v_i)=\vert E(S_{v_i}) \vert \geq 3,$ for $1 \leq i \leq k.$ 

\begin{lemma}[\cite{LT}]\label{Stardec}
If a tree $T$ has a star-decomposition and  $\vert E(T) \vert \equiv 1,2 \pmod 4 $, then $T$ is non-conservative.
\end{lemma}

Since a snowflake admits a star-decomposition, by Lemma \ref{Stardec} it follows that

\begin{remark}\label{ConservSnow}
If a snowflake of size $M$ is conservative then, $M \equiv 0,3 \ (mod \ 4)$.
\end{remark}

It can be seen that a $k$-$t$-conservative labeling of a snowflake of size $M$ represents a zero-sum $(M;k;t)$-Skolem system, and vice versa. So, in what follows, we will make no distinction between these structures.

\subsection{Even snowflakes}\label{Sect-4.1}

We say that a snowflake is  \textit{even} if every internal vertex, except possibly the center, has even degree.

\begin{remark}\label{etiquconsecu}
Let $S$ be a star of size $M$ with internal vertex $v$. Then, there exist Eulerian labelings $\phi, \phi_1$ and $\phi_2$ of $S$ such that $s(v)=0$ holds, where  
\begin{enumerate}
\item if $M\equiv \  0  \ (mod \ 4),$ then $\phi(E(S))=[1,M]$, 
\item if $M\equiv \  2 \ (mod \ 4),$ then $\phi_{1}(E(S))=[1,M-1]\cup\{M+1\}$ and $\phi_{2}(E(S))=\{1\}\cup[3,M+1].$ 
\end{enumerate} 
\end{remark}

\begin{lemma}\label{1copospar}
Let $\ell\in \{1,2\}$ and  $\mathbf{C}$ an even snowflake of size $M$ such that its center $z$ has degree two and $\{v_{1},v_{2}\}= V_{I}(\mathbf{C})\setminus\{z\}.$ Then, there exist Eulerian labelings $\phi, \phi_1$ and $\phi_2$ of $C$ such that $s(v_{1})=0=s(v_{2})$ and $s(z)=\ell$ hold, where

\begin{enumerate}
\item  if $M\equiv \  0  \ (mod \ 4),$ then $\phi(E(\mathbf{C}))=[1,M],$
\item if $M\equiv \  2 \ (mod \ 4),$ then $\phi_{1}(E(\mathbf{C}))=[1,M-1]\cup \{M+1\}$ and $\phi_{2}(E(\mathbf{C}))=\{1\}\cup [3,M+1].$
\end{enumerate}
\end{lemma}

\begin{proof}
Let $\ell\in \{1,2\}$, $\mathbf{C}$ an even snowflake of size $M$ such that its center $z$ has degree two and $G\cong S_{v_{1}}\cup S_{v_{2}}$ a galaxy obtained from $\mathbf{C}$ by splitting $z$ twice, where for $i\in\{1,2\},$ $S_{v_{i}}$ has size $M_{i}$, and $M=M_1+M_2$. Consider the following cases.
\\

\textit{Case 1.} $M \equiv 0 \ (mod \ 4).$

\textit{Case 1.1} $M_{i}\equiv 0 \ (mod \ 4).$

Set $\phi(E(S_{v_{1}})):= [1,M_{1}],$ and $\phi(E(S_{v_{2}})):= [M_{1}+1,M]$, where
$\phi(zv_{1})=M_{1}$ and $\phi(zv_{2})\in \{M_{1}+1,M_{1}+2\}$.

\textit{Case 1.2} $M_{i}\equiv 2 \ (mod \ 4).$

 Set $\phi(E(S_{v_{1}})):= [1,M_{1}-1]\cup\{M_{1}+1\},$ and $\phi(E(S_{v_{2}})):= \{M_{1}\}\cup [M_{1}+2,M]$, where $\phi(zv_{1})\in \{M_{1}-1,M_{1}-2\}$ and $\phi(zv_{2})=M_{1}$.
 
In either case, $\mathbf{C}$ admits an orientation such that $s(z)=\ell$ and by Remark \ref{etiquconsecu}  $s(v_{1})=0=s(v_{2}).$  
\\

\textit{Case 2.} $M \equiv 2 \ (mod \ 4).$ Assume $M_{1}\equiv 0 \ (mod \ 4)$  and $M_{2}\equiv 2 \ (mod \ 4).$
\\

Set $\phi_1(E(S_{v_{1}})):= [1,M_{1}]$ and $\phi_1(E(S_{v_{2}})):= [M_{1}+1,M_{1}+M_{2}-1]\cup\{M_{1}+M_{2}+1\}$, where $\phi_1(zv_{1})=M_{1}$ and $\phi_1(zv_{2})\in \{M_{1}+1,M_{1}+2\}$.

Set $\phi_2(E(S_{v_{1}})):= [M_{2}+2,M_{1}+M_{2}+1]$ and  $\phi_2(E(S_{v_{2}})):= \{1\}\cup [3,M_{2}+1]$, where $\phi_2(zv_{1})=M_{2}+2$ and $\phi_2(zv_{2})\in\{M_{2},M_{2}+1\}$.

In either case, $\mathbf{C}$ admits an orientation such that $s(z)=\ell$ and by Remark \ref{etiquconsecu}  $s(v_{1})=0=s(v_{2}).$ 

\end{proof}

\begin{figure}[H]
\begin{center}
\includegraphics[scale=.3]{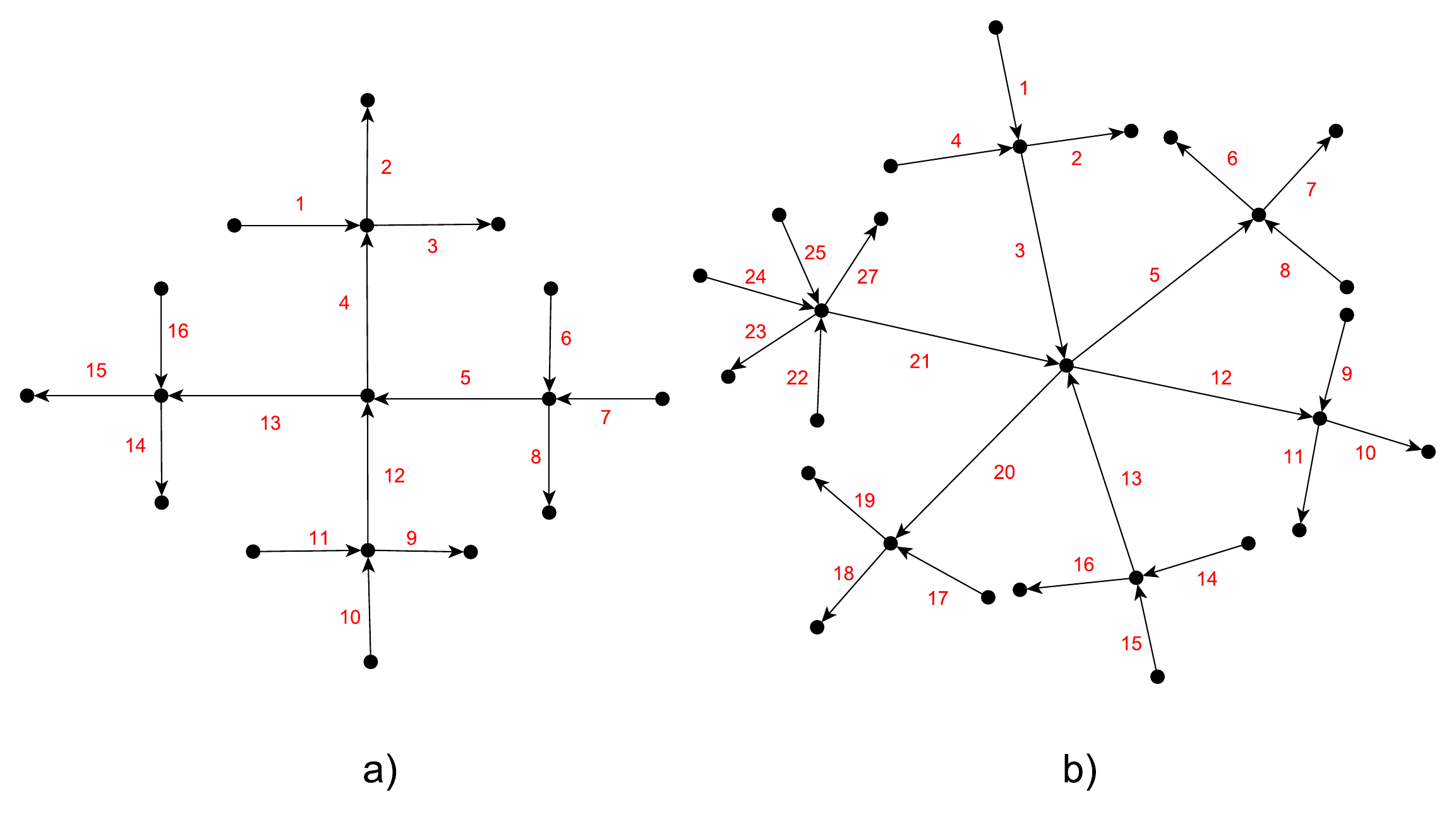} \caption{a) An Eulerian conservative labeling. b) An Eulerian near-conservative labeling.}
\label{ejemcopopar}
\end{center}
\end{figure}

\begin{lemma}\label{copospareseulerianoconser}
Let $\mathbf{C}$ be an even snowflake of size $M$ and $z$ its center such that $d_ \mathbf{C}(z)\equiv 0\  (mod \ 2).$  If $M\equiv 0$ $(mod \ 4),$ then  $\mathbf{C}$ admits an Eulerian conservative labeling. Otherwise $C$ admits an Eulerian near-conservative labeling. 
\end{lemma} 
\begin{proof}
Let $d_ \mathbf{C}(z)=n=2p$, $G$ be the galaxy obtained from $\mathbf{C}$ by splitting $z$ as many times as its degree, and $N_2= \vert \{ S \in G : \vert E(S) \vert \equiv 2 \pmod 4 \}\vert$. Notice that $\mathbf{C}=\mathbf{C}_1(z_1) \bullet \cdots \bullet \mathbf{C}_p(z_p)$, where $\mathbf{C}_i(z_i)$ is a snowflake of size $M_i$, centered at the vertex $z_i$ with $d_{C_i}(z_i)=2$. Consider the following cases.
%Then $\mathbf{C}$ results from the attachment of $n$  stars $S_1,\dots,S_n$ such that $\vert E(S_i) \vert =M_i$  and $M=\sum_{i=1}^n M_i$. 
\\

\textit{Case 1.} Let $M\equiv 0$ $(mod \ 4).$

Since $N_2$ is even, we can choose $\mathbf{C}_i$ such that for $1 \leq i \leq p$, $M_i \equiv 0 \pmod 4$. By Lemma \ref{1copospar}, we can obtain an Eulerian labeling $\phi'$ of $C$ such that if $p \equiv 0 \pmod 2$, then $s_{\overrightarrow{C_1},\phi'}(z_1)=1$, and  $s_{\overrightarrow{C_1,}\phi'}(z_1)=2$ otherwise. And for $2 \leq i \leq p$, $s_{\overrightarrow{C_i},\phi'}(z_i)=1$. By using properly remarks \ref{sumapar} and \ref{HS-HS}, we can obtain from $\phi'$ an Eulerian conservative labeling $\phi$ of $C$(see Figure \ref{ejemcopopar}a)).
\\

\textit{Case 2.} Let $M\equiv 2$ $(mod \ 4).$
\vspace{.2cm}

Since $N_2$ is odd, we can choose $\mathbf{C}_i$ such that for $1 \leq i \leq p-1$, $M_i \equiv 0 \pmod 4$ and $M_p  \equiv 2 \pmod 4$. The rest of the proof is analogous to that of Case 1 and will be omitted. (see Figure \ref{ejemcopopar}b)). 
\end{proof}

\begin{lemma}\label{coposparconserv}
Let $\mathbf{C}$ be an even snowflake of size $M$ and $z$ its center such that $d_ \mathbf{C}(z)\equiv 1\  (mod \ 2).$  If $M\equiv 0$ $(mod \ 4),$ then  $\mathbf{C}$ admits a  conservative labeling. Otherwise $C$ admits a near-conservative labeling.
\end{lemma} 

\begin{figure}[H]
\begin{center}
\includegraphics[scale=.3]{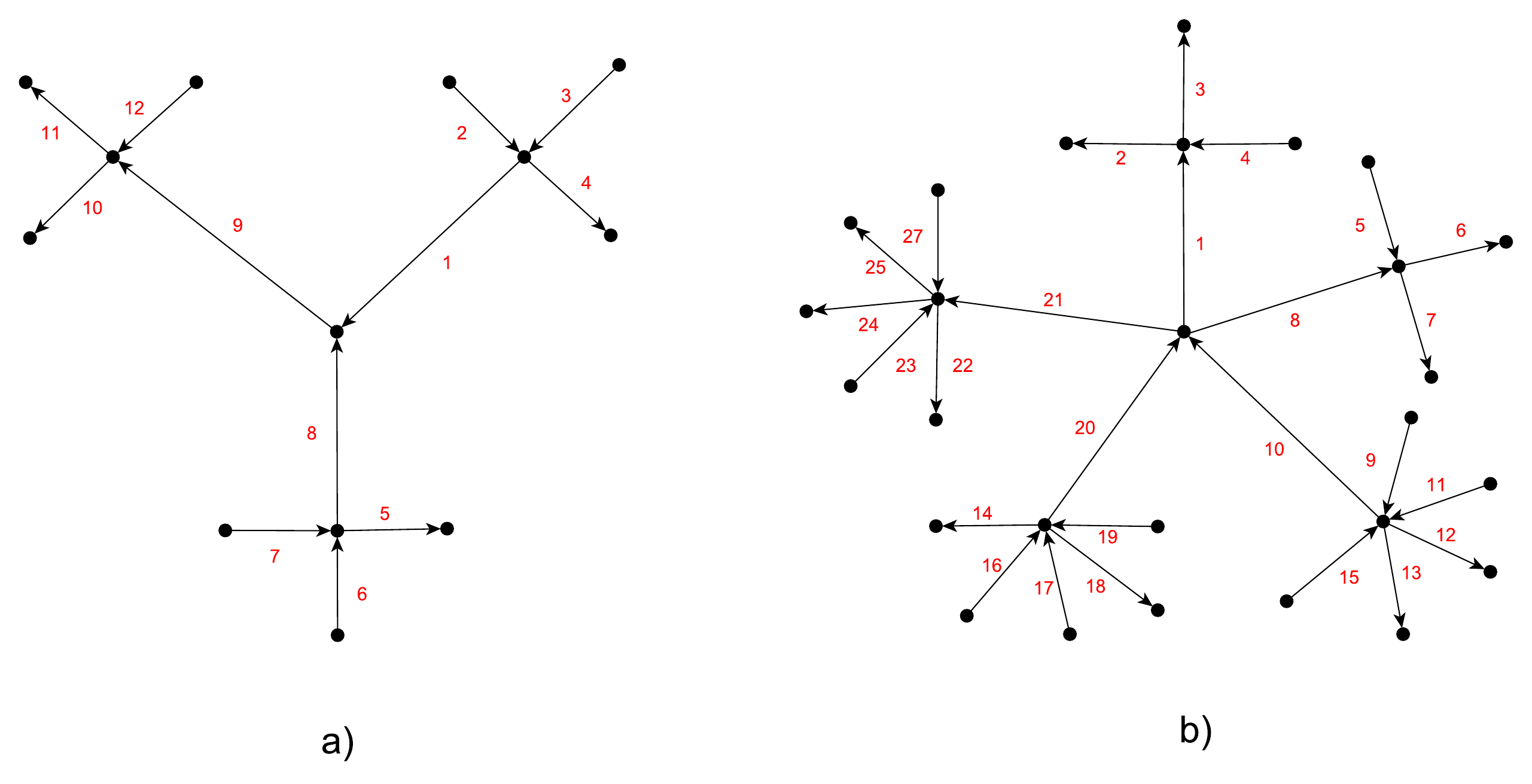} \caption{a) A conservative labeling. b) A near-conservative labeling.}
\label{ejemcopoimpar}
\end{center}
\end{figure}

\begin{proof}
Let $d_ \mathbf{C}(z)=n=2p+1$, $G$ be the galaxy obtained from $\mathbf{C}$ by splitting $z$ as many times as its degree, and $N_2= \vert \{ S \in G : \vert E(S) \vert \equiv 2 \pmod 4 \}\vert$. Notice that $\mathbf{C}=\mathbf{S}(z_0) \bullet \mathbf{C}_1(z_1) \bullet \cdots \bullet \mathbf{C}_p(z_p)$, where $S \in G$ with interior vertex $v$, $z_0$ is a leaf of $S$, and $\mathbf{C}_i(z_i)$ is a snowflake of size $M_i$, centered at the vertex $z_i$ with $d_{C_i}(z_i)=2$. Consider the following cases.
\\

\textit{Case 1.} Let $M\equiv 0$ $(mod \ 4).$
\vspace{.2cm}

Since $N_2$ is even, we can choose $\mathbf{C}_i$ such that for $1 \leq i \leq p$, $M_i \equiv 0 \pmod 4$ and $\vert E(S) \vert \equiv 0 \pmod 4$. By Remark \ref{etiquconsecu}, $S$ admits an Eulerian conservative labeling $\phi_S$, where $\phi_S(z_0v)=1$. Also, by Lemma \ref{1copospar}, we can obtain an Eulerian labeling $\phi'$ of $C$ such that if $p \equiv 0 \pmod 2$, then $s_{\overrightarrow{C_1},\phi'}(z_1)=2$, and  $s_{\overrightarrow{C_1,}\phi'}(z_1)=1$ otherwise. And for $2 \leq i \leq p$, $s_{\overrightarrow{C_i},\phi'}(z_i)=1$. By using properly remarks \ref{sumapar} and \ref{HS-HS}, we can obtain from $\phi_S$ and $\phi'$ an Eulerian conservative labeling $\phi$ of $C$(see Figure \ref{ejemcopoimpar}a)).
\\

\textit{Case 2.} Let $M\equiv 2$ $(mod \ 4).$
\vspace{.2cm}

Since $N_2$ is odd, we can choose $\mathbf{C}_i$ such that for $1 \leq i \leq p-1$, $M_i \equiv 0 \pmod 4$ and $M_p  \equiv 2 \pmod 4$. The rest of the proof is analogous to that of Case 1 and will be omitted. (see Figure \ref{ejemcopoimpar}b)). 
\end{proof}

The following theorem follows from Remark \ref{ConservSnow} and lemmas \ref{copospareseulerianoconser} and \ref{coposparconserv}. 
\begin{theorem}\label{princpares}
An even snowflake of size $M$ is conservative if $M\equiv 0 \ (mod \ 4),$ and it is near-conservative otherwise.
\end{theorem}

\subsection{Minimum odd snowflakes}\label{Sect-4.2}

We say that a snowflake is  \textit{odd} if every internal vertex, except possibly the center, has odd degree.

In this section we treat a $k$-$t$-conservative labeling of a snowflake of size $M$ as  a zero-sum $(M;k;t)$-Skolem system.

\begin{lemma}\label{15}
Let $k\geq 0,$ $n\geq 3,$ $M=3n$  and $t\in \{0,1\}.$ If $M\equiv t,3-t$ $(mod \ 4),$ then there exists a zero-sum $(3\cdot n;k;t)$-Skolem system of size $M$ and order $n.$ 
\end{lemma}
\begin{proof}
Let $t\in\{0,1\}$ and $S=\{D_{1},\ldots,D_{n}\}$ be the $(3\cdot n;k;t)$-Skolem system obtained by using the Skolem-O'Keefe Method. 
\\  

{\em Case 1.} Suppose that $n\equiv 0,3 \ (mod \ 4).$ %For $1\leq i\leq n,$ $k\geq 0$ and $t\in\{0,1\}$, let $S=\{D_{1},\ldots,D_{n}\}$ be a $(3\cdot n;k;t)$-Skolem system, we can assume that $D_{i}=\{(-i,-(n+a_{i}+k),n+b_{i}+k)\}.$
\vspace{.2cm}
 
Let $S'=abs(S)$. Notice that,  by construction, the numbers  $1,2,\ldots,n$ belong to distinct sets in $S'$. Also, by Remark \ref{star-cons}, the star of size $n$ is conservative. Hence, by using properly Remark \ref{HS-HS}, we can obtain a zero-sum  $(3\cdot n;k;t)$-Skolem system $S_0$ from $S$, where $Z_0(S_0)$ satisfies $abs(Z_0(S_0))= \{1,2,\ldots,n \}$.
\\

{\em Case 2.} Suppose that  $n\equiv 1 \ (mod \ 4).$ 
\vspace{.2cm}

Let $n=4s+1.$ For $s=1$ and $s=2$, see  Figure \ref{ct9}. By Lemma \ref{ordenR1}, $S$ admits an $R$-order say $S_R=(D_{1},\ldots,D_{n})$.

\begin{figure}[H]
\begin{center}
\includegraphics[scale=.4]{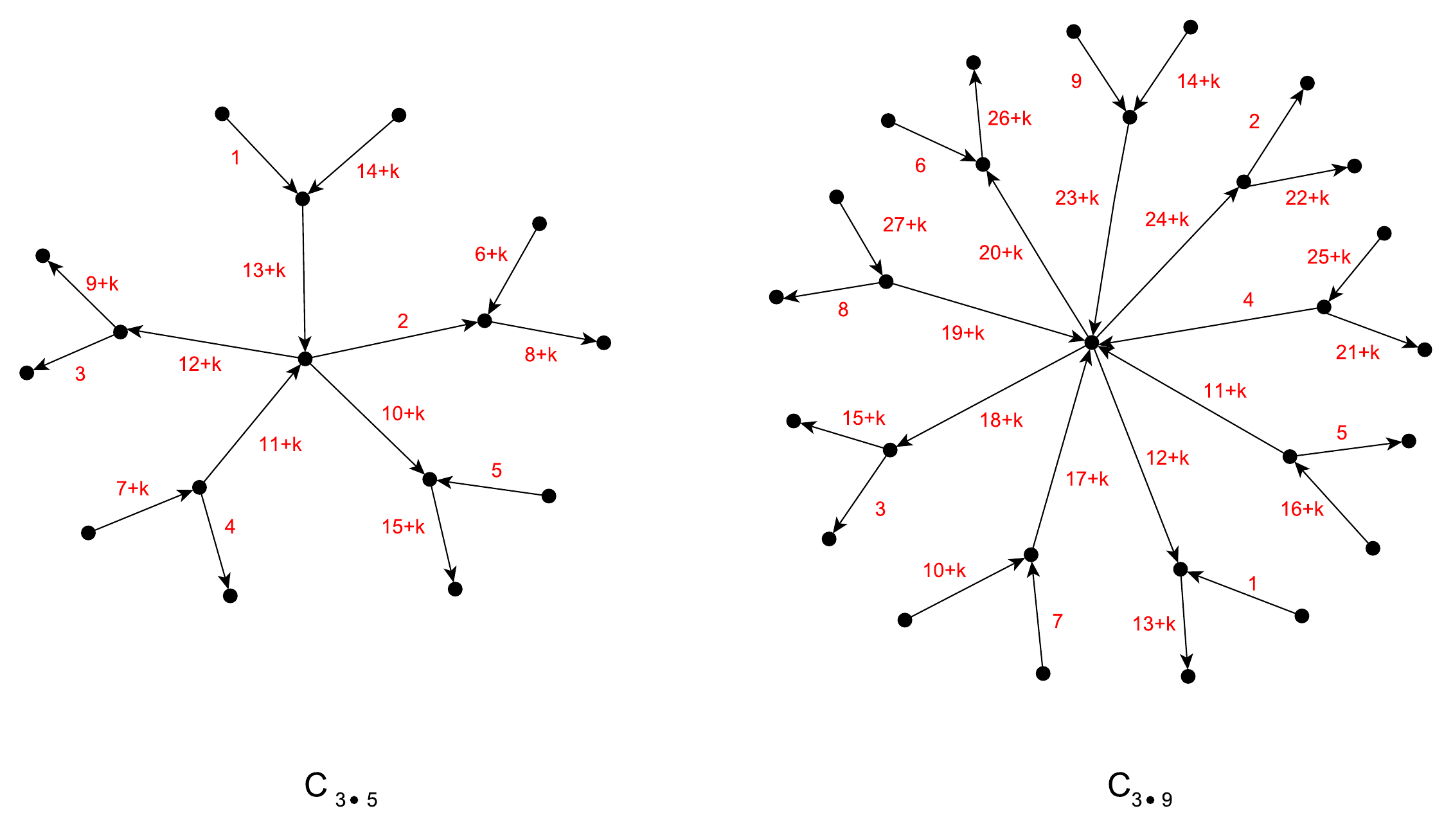} \caption{ Zero-sum Skolem systems with $n=5$ (left) and  $n=9$ (right) (viewed as snowflakes).}
\label{ct9}
\end{center}
\end{figure}

Let $s\geq 3$.  
Notice that, by Remark  \ref{ordenRcong3}, $D_{2s+3}=\{(-4,-(n+6s+k),n+6s+4+k\}$. Moreover, for all $1\leq i\leq 2s+2$ and $2s+4\leq i\leq 4s+1$ there exist $d_{i}\in abs(D_{i})$ and $d_{i+1}\in abs(D_{i+1})$ such that $|d_{i}-d_{i+1}|=1$. So, we obtain that 
\\
$4+\sum\limits_{i=1}^{\frac{2s+2}{2}} |d_{2i-1}-d_{2i}|+\sum\limits_{i=\frac{2s+5}{2}}^{\frac{4s+1}{2}}|d_{2i-1}-d_{2i}|=2s+4\equiv 0 \ (mod \ 2)$.   
\\
Hence, by Remark \ref{sumapar} and by using properly Remark \ref{HS-HS}, we can obtain a zero-sum  $(3\cdot n;k;t)$-Skolem system $S_0$ from $S$, where $Z_0(S_0)$ satisfies $abs(Z_0(S_0))= \{4 \} \cup \left( \bigcup_{i \neq 2s+3} d_i \right)$.
\\

\textit{Case 3.} Suppose that $n\equiv 2 \ (mod \ 4).$ 
\vspace{.2cm}

Let $n=4s+2.$ For $s=1$, see Figure \ref{ct(6)}. 

\begin{figure}[H]
\begin{center}
\includegraphics[scale=.4]{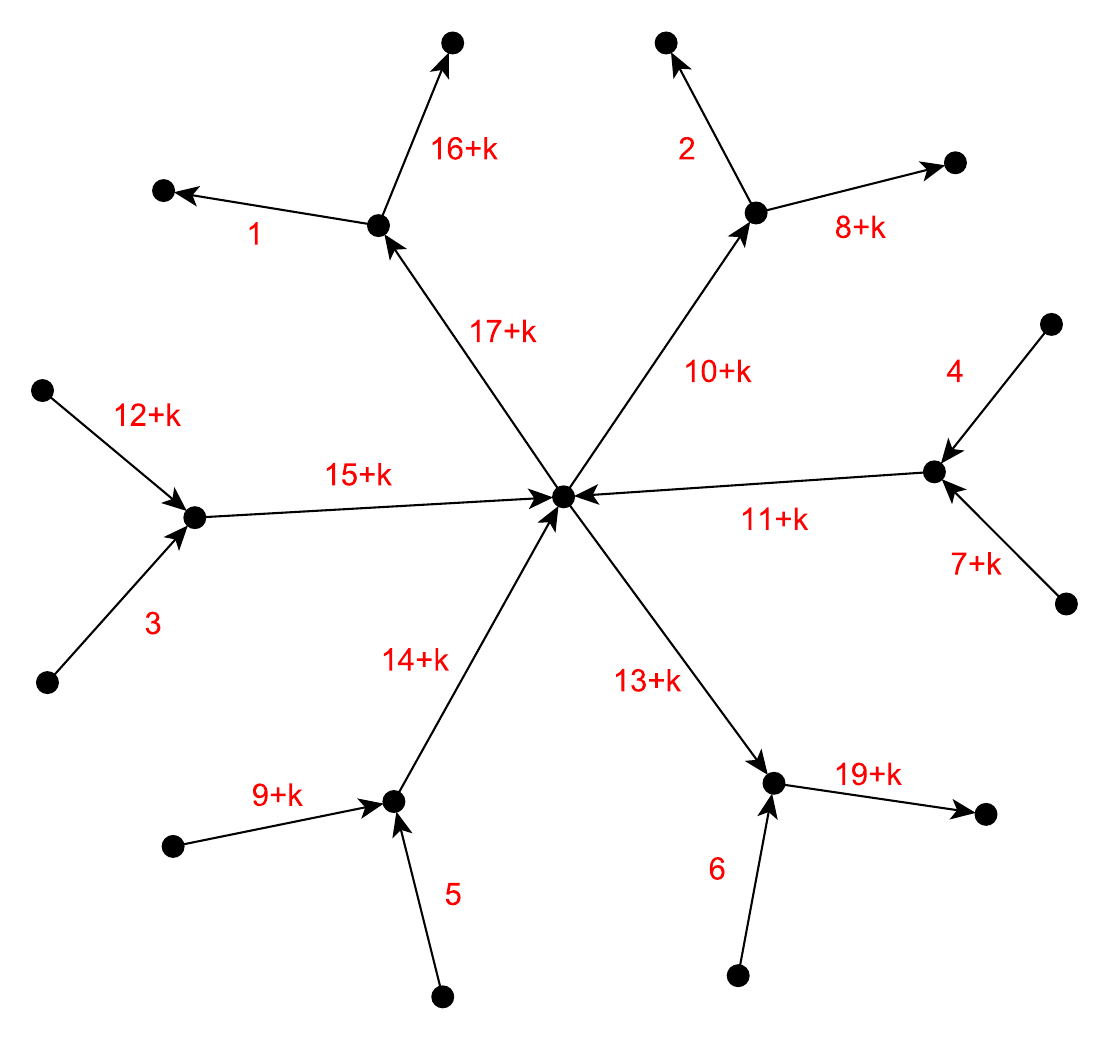} \caption{A zero-sum $(3\cdot 6;k;1)$-Skolem system (viewed as a snowflake).}\label{ct(6)}
\end{center}
\end{figure}

For $s\geq 2,$ we consider two cases.
\\

\textit{Case 3.1.}  Suppose $s$ is even. 
\vspace{.2cm}
By Lemma \ref{ordenR2}, for some $q$  odd, there exist $d\in abs(D_{q})$ and $d'\in abs(D_{q+1})$  such that  $|d-d'|=2$. In fact, for $q=3s+1$ take $d=n+7s+2+k$ and $d'=n+7s+4+k$ as needed. Moreover, for $S\setminus \{D_{3s+1},D_{3s+2}\}$, for all $1\leq i\leq 3s$ and $3s+3\leq i\leq 4s+2$ there exist $d_{i}\in D_{i}$ and $d_{i+1}\in D_{i+1}$ such that $|d_{i}-d_{i+1}|=1$. 
\\
Thus, proceeding as in Case 2, we can obtain a zero-sum  $(3\cdot n;k;t)$-Skolem system $S_0$ from $S$, where $Z_0(S_0)$ satisfies $abs(Z_0(S_0))= \{d,d' \} \cup \left( \bigcup_{i \notin \{3s+1, 3s+2 \}} d_i \right)$.
\\

\textit{Case 3.2.}  Suppose $s$ is odd. The proof of this case is similar to that of Case 3.1. Just take into account that $q=3s$, $d=n+7s+1+k$ and $d'=n+7s+3+k$.
\end{proof}

\begin{lemma}\label{16}
Let $k\geq 0,$ $n=n_{1}+n_{2}\geq 3,$ $M=3n_ {1}+5n_{2}$  and $t\in \{0,1\}.$ If $M\equiv t,3-t$ $(mod \ 4),$ then there exists a zero-sum $(3\cdot n_{1},5\cdot n_{2};k;t)$-Skolem system of size $M$ and order $n.$ 
\end{lemma}
\begin{proof}
Let $k\geq 0,$ $t\in\{0,1\}$ and $S$ be the $(3\cdot n_{1},5\cdot n_{2};k;t)$-Skolem  system  obtained in the proof of Lemma \ref{3,5hefftersystem}. Notice that the set of distinguish elements obtained for the zero-sum $(3\cdot n;k;t)$-Skolem  system constructed in the proof of  Lemma \ref{15} can be taking as a set of distinguish element in $S$, in order to construct a zero-sum $(3\cdot n_{1},5\cdot n_{2};k;t)$-Skolem system.
\end{proof}

%Figure \ref{C(3,3,3,5)} is an example of $\phi:E(C_{3(3),5(1)})\rightarrow [1,4]\cup[5+k,13+k]\cup\{15+k\},$ a $k$-$1$-conservative labeling  obtained from Lemma \ref{16}, where $S=\{(1,7+k,8+k),(2,11+k,13+k),(3,9+k,12+k),(4,10+k,15+k,5+k,6+k)\}$ and $r(S)=4.$ 
\begin{figure}[H]
\begin{center}
\includegraphics[scale=.5]{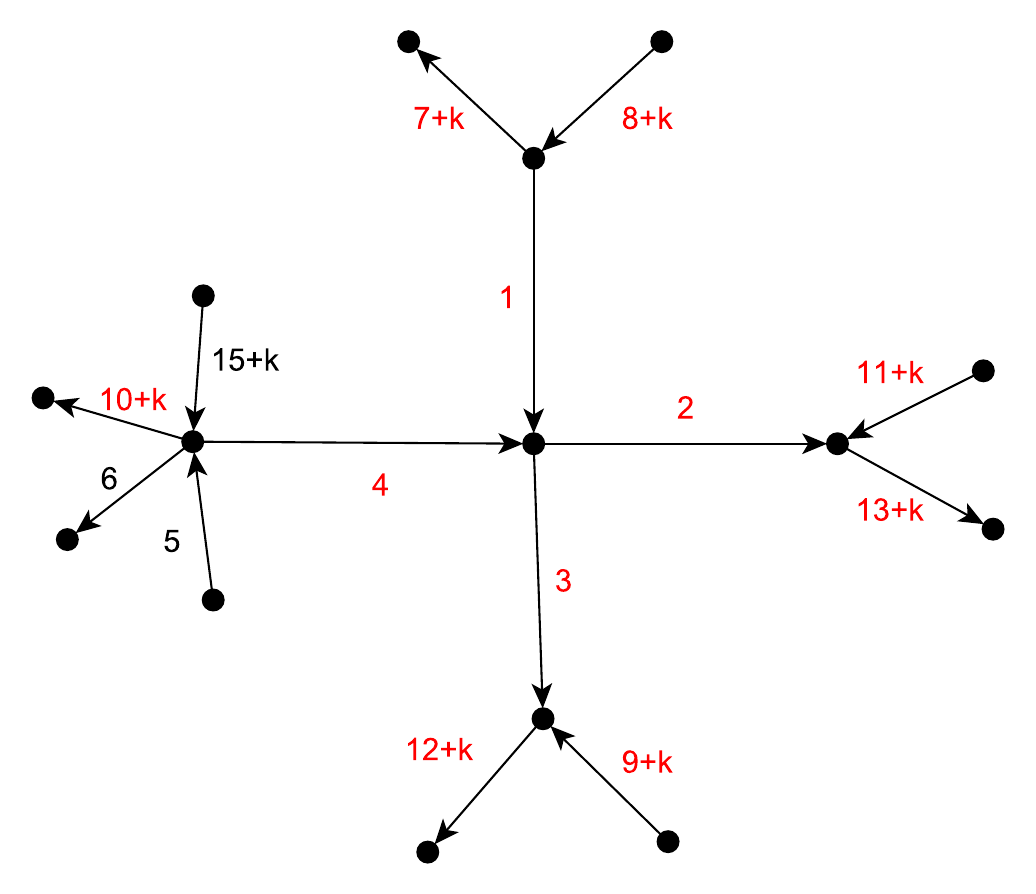}\caption{ Zero $(3\cdot 3,5\cdot 1,;k;1)$-Skolem system.}\label{C(3,3,3,5)}
\end{center}
\end{figure}

\subsection{Attaching minimum odd snowflakes with even snowflakes}\label{Sect-4.3}

\begin{lemma}\label{2copospar}
Let $x\in\mathbb{N},$  $\ell\in \{x-1,x,x+1\},$  $\mathbf{C}$ an even snowflake of size $M$ such that its center $z$ has odd degree $n \geq 3$  and $\{v_i\}_{1}^{n}= V_{I}(\mathbf{C})\setminus\{z\}.$ Then there exist labelings $\phi_0, \phi_1, \phi_2$ and $\phi_3$, and orientations $\overrightarrow{\mathbf{C}}_0$, $\overrightarrow{\mathbf{C}}_1$, $\overrightarrow{\mathbf{C}}_2$ and $\overrightarrow{\mathbf{C}}_3$ of $\mathbf{C}$ such that $s_{\overrightarrow{\mathbf{C}}_i,\phi_i}(v_{i})=0$ and $s_{\overrightarrow{\mathbf{C}}_i,\phi_i}(z)=\ell$ hold, where
\begin{enumerate}
\item  if $M\equiv \  0  \ (mod \ 4),$ then $\phi_0(E(\mathbf{C}))=[x,M+x-1],$
\item  if $M\equiv \  0  \ (mod \ 4)$ and $\mathbf{C}=\mathbf{C}'(z)\bullet S_4(u)$, where $S_4$ is a star of size four and $u$ is a leaf of $S_4$, then $\phi_1(E(\mathbf{C}))= \phi_1(E(\mathbf{C}')) \cup \phi_1(E(S_4))=[r+1,M+r-4]\cup\{x\}\cup \{x+2,x+3,x+5\}$, where $0 < r < x-M+4$.
\item if $M\equiv \  2 \ (mod \ 4),$ then $\phi_{2}(E(\mathbf{C}))=[x,M+x-2]\cup \{M+x\}$ and $\phi_{3}(E(\mathbf{C}))=\{x\}\cup [x+2,M+x].$
\end{enumerate} 
\end{lemma}

\begin{proof}
The proof of this lemma is easy due to Remark \ref{etiquconsecu} and  Lemma \ref{1copospar}.  
\end{proof}

\begin{lemma}\label{1coposimpar}
Let $k\geq 0,$ $t\in\{0,1\},$ $\mathbf{C}=C_{3\cdot n_{1},5\cdot n_{2}}$ an odd snowflake of size $M\equiv t,3-t$ $(mod \ 4)$ such that its center $z$ has degree $n\geq 1$ and $\{v_{i}\}_{i=1}^{n}= V_{I}(\mathbf{C})\setminus\{z\}$.  Then, there exists a labeling  $\phi:E(\mathbf{C})\rightarrow [1,n]\cup[n+1+k,M+k-1]\cup\{M+k+t\}$  and an orientation of $\mathbf{C}$  such that   $s_{\overrightarrow{\mathbf{C}}}(v_{i})=0$ and $s_{\overrightarrow{\mathbf{C}}}(z)\in\{1,2\}.$  

\end{lemma}

\begin{proof}
By Lemma \ref{3,5hefftersystem}, there exists a  $(3\cdot n_{1},5\cdot n_{2};k;t)$-Skolem system of size $M\equiv t,3-t$ $(mod \ 4)$, $S=\{D_{1},\ldots,D_{n}\}$. Let $N=\{1,\ldots,n \}$. Notice that the numbers in $N$ belong to distinct $abs(D_{i})$. Hence $S$ induce a labeling $\phi$ of $\mathbf{C}=C_{3\cdot n_{1},5\cdot n_{2}}$, where for $1\leq i\leq n,$  $\phi(zv_{i})=i$ and $s(v_{i})=0.$ 

On the other hand, it is clear that $x \in N$ admits a sign, $sign(x)$, such that $s_N=\sum_{i=1}^n sign(i)=2$ if $n\equiv 0,3 \  (mod \ 4),$ and $s_N=1$ otherwise. Hence, by using Remark \ref{HS-HS}, $\mathbf{C}$ admits an orientation $\overrightarrow{\mathbf{C}}$ such that $s_{\overrightarrow{\mathbf{C}}}(v_{i})=0$ and $s_{\overrightarrow{\mathbf{C}}}(z)=s_N\in\{1,2\}.$  
\end{proof}

\begin{lemma}\label{2coposimpar}
Let $k\geq 0,$ $t\in\{0,1\},$   $\mathbf{C}=C_{3\cdot n_{1},5\cdot n_{2}}$ an odd snowflake of size $M\equiv t,3-t$ $(mod \ 4)$ such that its center $z$ has degree $n\geq 1$ and $\{v_{i}\}_{i=1}^{n}= V_{I}(\mathbf{C})\setminus\{z\}.$ Then, there exists a bijective labeling  $\phi:E(\mathbf{C})\rightarrow [1,n]\cup[n+1+k,M+k-1]\cup\{M+k+t\}$ and an orientation of $\mathbf{C}$  such that   $s_{\overrightarrow{\mathbf{C}}}(v_{i})=0$ and $s_{\overrightarrow{\mathbf{C}}}(z)\in  \{M+k+t,M+k+t+1,M+k+t+2\}.$   

\end{lemma}

\begin{proof}
The proof is analogous to the proof of Lemma \ref{1coposimpar} but considering that:
\\
for $1\leq i\leq n-2,$ 
$$\phi(zv_{i})=i.$$
$$\phi(zv_{n-1})\in\{n,n-1\}.$$
$$\phi(zv_{n})\in\{M,M+1\}.$$

Thus, we can obtain that,
\begin{itemize}
\item if $n\equiv 0,1 \  (mod \ 4)$ and $\{n,M+k+t\}\in D_{n}$, then $s_{\overrightarrow{\mathbf{C}}}(z)=M+k-t+2,$
\item if $n\equiv 0 \  (mod \ 4)$ and $\{n-1,M+k+t\}\in D_{n}$, then $s_{\overrightarrow{\mathbf{C}}}(z)=M+k+t+1,$
\item if $n\equiv 2,3 \  (mod \ 4)$ and $\{n,M+k+t\}\in D_{n}$, then $s_{\overrightarrow{\mathbf{C}}}(z)=M+k+t+1,$
\item if $n\equiv 2 \  (mod \ 4)$ and $\{n-1,M+k-t\}\in D_{n}$, then $s_{\overrightarrow{\mathbf{C}}}(z)=M+k-t+2,$
\item if $n=1$ and $\{1,M+k+t\}\in D_{1}$, then $s_{\overrightarrow{\mathbf{C}}}(z)=M+k+t,$
\item if $n=2$ and $\{2,M+k+t\}\in D_{2}$, then $s_{\overrightarrow{\mathbf{C}}}(z)=M+k+t+1,$
\item if $n=2$ and $\{1,M+k+t\}\in D_{2}$, then $s_{\overrightarrow{\mathbf{C}}}(z)=M+k+t+2.$
\end{itemize}
\end{proof}

\begin{lemma}\label{attachoddeven2}
Let $k\geq 0,$ $t\in \{0,1\}$ and $\mathbf{C}=C_{3\cdot n_{1},5\cdot n_{2},6\cdot n_{3},4\cdot n_{4}}$ a snowflake of size  $M=3n_ {1}+5n_{2}+6n_{3}+4n_{4}$ such that $n_{1}+n_{2}+n_{3}+n_{4}\geq 3$, $n_{1}+n_{2}\geq 1$ and $ n_{3}+n_{4}\geq 1.$ If $M\equiv t,3-t$ $(mod \ 4),$ then $C$ admits a $k$-$t$-conservative labeling. 
\end{lemma}

\begin{proof}
Let $t'\in \{0,1\},$ $\mathbf{C}=C_{3\cdot n_{1},5\cdot n_{2},6\cdot n_{3},4\cdot n_{4}}$ of size $M,$ $\mathbf{C_{1}}=C_{3\cdot n_{1},5\cdot n_{2}}$ of size $M_{1}\equiv t',3-t'$ with $n=n_{1}+n_{2},$ $z_{1}$ its center and $\{v_{i}\}_{i=1}^{n}= V_{I}(\mathbf{C_{1}})\setminus\{z_{1}\}$,  $\mathbf{C_2}=C_{6n_{3},4n_{4}}$ of size $M_{2}$ with $n'=n_{3}+n_{4},$ $z_{2}$ its center and $\{w_{i}\}_{i=1}^{n'}= V_{I}(\mathbf{C_{2}})\setminus\{z_{2}\}.$ Notice that  $\mathbf{C}\cong \mathbf{C_1}(z_{1})\bullet\mathbf{C_2}(z_{2})$ where $M=M_{1}+M_{2}$ and $z=z_{1}=z_{2}$ is the center of $\mathbf{C}.$  Now, we consider the following cases.
\\

\textit{Case 1.} $ n_{3}+n_{4}\equiv 0 \ (mod \  2).$  
\\

\textit{Case 1.1.} Suppose that $n_3$ is odd.
\vspace{.1cm}

By Lemma \ref{1coposimpar},  there is  a labeling $\phi_{1}:E(\mathbf{C_{1}})\rightarrow[1,n]\cup[n+1+k,M_{1}+k-1]\cup\{M_{1}+k+t'\}$ and an orientation  of $\mathbf{C_{1}}$ such that  $s_{\overrightarrow{\mathbf{C_{1}}}}(v_{i})=0$ and  $s_{\overrightarrow{\mathbf{C_{1}}}}(z_{1})=\ell\in\{1,2\}.$
\\

By applying properly  Lemma \ref{1copospar} as many times as needed we obtain a bijective labeling  $\phi_{2}:E(\mathbf{C_{2}})\rightarrow [M_{1}+k+1,M+k-1]\cup\{M+k+1\}$ if $t'=0,$  or $\phi_{2}:E(\mathbf{C_{2}})\rightarrow \{M_{1}+k\}\cup[M_{1}+k+2,M+k]$  if $t'=1.$ In both cases there is an orientation of $\mathbf{C_{2}}$  such that  $s_{\overrightarrow{\mathbf{C_{2}}}}(w_{i})=0$ and $s_{\overrightarrow{\mathbf{C_{2}}}}(z_{2})=\ell$.
\\
Thus, by reversing the arcs in  $\mathbf{C_{2}}$  we have $s_{\overrightarrow{\mathbf{C_{2}}}}(z_{2})=-\ell.$ Hence, $s_{\overrightarrow{\mathbf{C}}}(z)=0.$   
\\

Therefore,  we can construct  a  $k$-$t$-conservative labeling $\phi$ of $\mathbf{C},$ where $\phi_{E(C_1)}=\phi_1$ and  $\phi_{E(C_2)}=\phi_2$, and $t=1-t'.$
\\

\textit{Case 1.2.} Suppose that $n_3$ is even.
\vspace{.1cm}

The proof is analogous to that of Case 1.1 but now consider that, $\phi_{1}:E(\mathbf{C_{1}})\rightarrow[1,n]\cup[M_{2}+n+1+k,M+k-1]\cup\{M+k+t\}$ and  $\phi_{2}:E(\mathbf{C_{2}})\rightarrow [n+1,M_{2}+n].$
\\
Therefore, $\phi_{1}$ and $\phi_{2}$ determine a $k$-$t$-conservative labeling of $\mathbf{C},$ where $t=t'.$ 
\\

\textit{Case 2.} $ n_{3}+n_{4}\equiv 1 \ (mod \  2).$
\vspace{.1cm}

\textit{Case 2.1.} $ n_{3}+n_{4} \geq 3$.

By Lemma \ref{2coposimpar}, there is  a labeling $\phi_{1}:E(\mathbf{C_{1}})\rightarrow[1,n]\cup[n+1+k,M_{1}+k-1]\cup\{M_{1}+k+t'\}$ and an orientation  of $\mathbf{C_{1}}$ such that  $s_{\overrightarrow{\mathbf{C_{1}}}}(v_{i})=0$ and  $s_{\overrightarrow{\mathbf{C_{1}}}}(z_{1})=\ell\in \{M+k+t',M+k+t'+1,M+k+t'+2\}$.

\textit{Case 2.1.1} Suppose that $n_3$ is odd. In this case $M_2 \equiv 2 \pmod 4$.
\vspace{.1cm}
The rest of the proof of this case is analogous to that of Case 1.1 but now consider the labelings described in Lemma \ref{2copospar}.3. 
\\
%Therefore,  $\phi_{1}$ and $\phi_{2}$ determine a $k$-$t$-conservative labeling of $\mathbf{C},$ where $t=1-t'.$

\textit{Case 2.1.2.} Suppose that $n_3$ is even. In this case $M_2 \equiv 0 \pmod 4$.

\vspace{.1cm}

By Lemma \ref{2copospar}.1, when $t'=0$  there is labeling  $\phi_2:E(\mathbf{C_{2}})\rightarrow [M_{1}+k+1,M+k]$ and an orientation of $C$ such that  $s_{\overrightarrow{\mathbf{C_{2}}}}(w_{i})=0$ and $s_{\mathbf{C_{2}}}(z_{2})=\ell.$

When $t'=1$, consider that $\mathbf{C_2}=\mathbf{C_2'}(z_2)\bullet\mathbf{S}_v(w)$, where $\mathbf{S}_v(w)$ is a star of size four centered at $v$ and $w$ is a leaf of $\mathbf{S}_v(w)$. In this case there exists $\phi_{1}:E(\mathbf{C_{1}})\rightarrow[1,n]\cup[M_2+n+k-3,M+k-5]\cup\{M+k-3\}$. By Lemma \ref{2copospar}.2, there is labeling   $\phi_{2}:E(\mathbf{C_{2}})\rightarrow [n+1,M_{2}+n-4]\cup\{M+k-4,M+k-2,M+k-1,M+k+1\}.$ Here $\phi_2(E(\mathbf{S}_v(w)))=\{M+k-4,M+k-2,M+k-1,M+k+1\}$.
\\ 
In both cases  there is an orientation  of $\mathbf{C_{2}}$ such that  $s_{\overrightarrow{\mathbf{C_{2}}}}(w_{i})=0$ and $s_{\mathbf{C_{2}}}(z_{2})=\ell.$
\\
Thus, by reversing the arcs in  $\mathbf{C_{2}}$  we can obtain that $s_{\overrightarrow{\mathbf{C_{2}}}}(z_{2})=-\ell.$ Hence, $s_{\overrightarrow{\mathbf{C}}}(z)=0.$    
\\
   
Therefore,  using $\phi_{1}$ and $\phi_{2}$, we can construct a $k$-$t$-conservative labeling of $\mathbf{C},$ where $t=t'.$

\textit{Case 2.2.} $ n_{3}+n_{4} =1$.
\\
This case is treated in the same way as Case 2.1. Just take into account that $n_1+n_2 \geq 2$ and $s_{\overrightarrow{\mathbf{C_{1}}}}(z_{1})=\ell\in \{M+k+t'+1,M+k+t'+2\}$.

\end{proof}

\begin{lemma}\label{PRINC}
Let $t \in \{0,1\}$. A snowflake of size $M \equiv t,3-t$ $(mod \ 4)$ admits a $t$-conservative labeling.
\end{lemma} 

\begin{proof}
Let $\mathcal{C}$ be a snowflake of size $M.$ If $\mathcal{C}$  is an even  snowflake, then by Theorem \ref{princpares} the result holds. In other case, $\mathcal{C}$ results from a proper attachment of a snowflake  $\mathcal{C}'$ and $p$ stars $S_1,\ldots,S_p$, where
\begin{itemize}
\item $G=\bigcup_{i=1}^pS_i$ is a (possibly empty) galaxy such that $\vert E(S_i) \vert \equiv 0 \pmod 4$, for $1 \leq i \leq p$,
\item $\mathcal{C}'$ is a minimum snowflake with at least one vertex, distinct of its center, of odd degree.
\end{itemize}
Notice that, in order to obtain $\mathcal{C}$, each center of $S_i$ most be identified with a vertex of $\mathcal{C}'$, distinct of its center, in the aforth mentioned attachment.

Let $M_G= \sum_{i=1}^p \vert S_i \vert$ and $M'=\vert E(\mathcal{C}') \vert$. If $\mathcal{C}'$ is an odd snowflake, by Lemma \ref{16}, $\mathcal{C}'$ admits a $M_G$-$t$-conservative labeling $\phi'$, where $\phi'(E(\mathcal{C}')=[1,r_1]\cup[r_1+1+M_G,M'+M_G-1]\cup\{ M'+M_G+t\}$. Set $r=r_1$.

Otherwise, $\mathcal{C}'=C_{3\cdot n_{1},5\cdot n_{2},6\cdot n_{3},4\cdot n_{4}}$, where $n_{1}+n_{2}+n_{3}+n_{4}\geq 3$, $n_{1}+n_{2}\geq 1$ and $ n_{3}+n_{4}\geq 1$. By Lemma \ref{attachoddeven2}, $\mathcal{C}'$ admits a $M_G$-$t$-conservative labeling $\phi'$, where $\phi'(E(\mathcal{C}')=[1,r_2]\cup[r_2+1+M_G,M'+M_G-1]\cup\{ M'+M_G+t\}$. Set $r=r_2$.

On the other hand, by applying properly Remark \ref{etiquconsecu}, we can construct an Eulerian conservative labeling $\phi$ of the galaxy $G$. Let $\phi_{r}$ defined as $\phi_{r}(e):=\phi(e)+r$, notice that $\phi_r(E(G))=[r+1,r+M_G]$. By Remark \ref{fuerteconser}, $s_{\overrightarrow{G},\phi_r}(v)=0$ for every $v \in V_I(G)$.

Therefore, there exists a $t$-conservative labeling $\theta$ of $\mathcal{C}$ such that $\theta_{\mathcal{C}'}=\phi'$ and $\theta_{G}=\phi$.

\end{proof}

Taking into account lemmas \ref{ConservSnow} and \ref{PRINC} it is straightforward to prove the main result of this section.

\begin{theorem}\label{Snowconservative}
A snowflake of size $M$ is conservative if $M\equiv 0,3$ $(mod \ 4)$, and it is near-conservative otherwise.
\end{theorem}

In view of the present work a natural question arises: is graceful every two nested cycles graph of size $M \equiv 0,3 \pmod 4$, and is it near-graceful otherwise?


\begin{thebibliography}{XX}

\bibitem[Acharya & Gill, (1981)]{AG}
B. D. Acharya and M. K. Gill, 
\textit{On the index of gracefulness of a graph and the gracefulness of two-dimensional square lattice graphs,}
\newblock{Indian J. Math.,}
\newblock{ 23, 81-94, 1981.} 

\bibitem[Alkasasbeh & Dyer, (2021)]{AD}
A. H. Alkasasbeh and D. Dyer,
\textit{Graceful Labellings of Various Cyclic Snakes},
https://arxiv.org/pdf/2012.10341, 2021.

\bibitem[Archdeacon et al. (2017)]{ABD}
D. Archdeacon, T. Boothby and J. H. Dinitz,
\textit{Tight Heffter Arrays Exist for all Possible Values,}
\newblock {J. Combin. Des.}
\newblock {25 (1),  5-35, 2017}.

\bibitem[Archdeacon et al. (2015)]{Arch4} 
D.~Archdeacon, J.H Dinitz, D.M. Donovan and E. Yazici, 
\textit{Square integer Heffter arrays with empty cells}, 
\emph{Des. Codes Cryptogr.} 77, no. 2-3, 409-426, 2015.


\bibitem[Bange et al. (1980)]{BBS}
D. Bange,  A. Barkauskas and  P. Slater,
\textit{Conservative Graphs},
\newblock {Journal of Graph Theory},  
\newblock {Vol. 4},
\newblock {81-91,}
\newblock {1980.}

\bibitem[Barrientos, (2001)]{Ba0}
C. Barrientos, 
\textit{Graceful labelings of cyclic snakes,}
\newblock {Ars Combin.,}
\newblock {60,}
\newblock {85-96,}
\newblock {2001.}

\bibitem[Barrientos, (2015)]{Ba1}
C. Barrientos, 
\textit{On graceful chain graphs,}
\newblock {Util. Math.,}
\newblock {78},
\newblock{55-64,}
\newblock {2009.}

\bibitem[Delorme et al. (1980)]{DKMTT}
C. Delorme, K.M. Koh, M. Maheo, H.K. Teo and H. Thuillier,
\textit{Cycles with a chord are graceful,}
\newblock {J. Graph Theory 4},  
\newblock {409-415},
\newblock {1980.}

\bibitem[Dinitz & Wanless (2017)]{DW}
J.H. Dinitz and I.M. Wanless,
\textit{The existence of square integer Heffter arrays,}
\newblock{Ars. Math. Contemp.,}
\newblock{13},
\newblock{81-93},
\newblock{2017.}

\bibitem[Elumalai et al. (2015)]{EE}
A. Elumalai and A. Anand Ephremnath,
\textit{Gracefulness of a cycle with zigzag chords, }
\newblock {International Journal of Pure and Applied Mathematics,}
\newblock {2015.}

\bibitem[Frucht & Gallian, (1988)]{FG}
R. Frucht and J. A. Gallian, 
\textit{Labeling prisms,}
 \newblock {Ars Combin., 26, 69-82, 1988.}

\bibitem[Gallian, (1989)]{G}
J. A. Gallian, 
\textit{Labeling prisms and prism related graphs}, 
\newblock{Congr. Numer.,}
 \newblock{59, 89-100, 1989.}
 
\bibitem[Goldfeder & Tey, (2018)]{GT}
I. Goldfeder and J. Tey,
\textit{A note on conservative galaxies, Skolem systems, cyclic cycle decompositions, and Heffter arrays},
\newblock {Discrete Mathematics 341,}  
\newblock {2519-2528,}
\newblock {2018.}

\bibitem[Gnanajothi, 1991]{Gn}
R. B. Gnanajothi, 
\textit{Topics in Graph Theory,}
 \newblock {Ph. D. Thesis, Madurai Kamaraj Uni- versity,}
 \newblock {1991}.

\bibitem[Koh & Punnim, (1982)]{KP}
K.M. Koh and N. Punnim,
\textit{On graceful graphs: cycles with 3-consecutive chords,} 
\newblock {Bull. Malaysian Math. Soc.},
\newblock {5, 49-63, 1982.}

\bibitem[Licona & Tey, (2024)]{LT}
M. Licona and J. Tey,
\textit{Conservative trees,} 
\newblock {Discrete Mathematics},
\newblock {347, Issue 3, 113854, 2024.}

\bibitem[Liu, (1989)]{L}
Y. Liu,
\textit{On Bodendiek's conjecture for graceful graphs,} 
\newblock {Chinese Quart. J. Math.,}
\newblock {4 , suppl., 67-73, 1989.}

\bibitem[Ma, (1988)]{Ma}
X. Ma,
\textit{A graceful numbering of a class of graphs,} 
\newblock {J. Math. Res. and Exposition},
\newblock {215-216, 1988.}

\bibitem[Ma et al. (1990)]{MLL}
X. Ma, Y.Liu and W. Liu,
\textit{Graceful graphs: cycles with $(t-1)$ chords,} 
\newblock {Math. Appl.,}
\newblock {9 , suppl., 6-8, 1990.}

\bibitem[Maheo, (1980)]{M}
M. Maheo, 
\textit{Strongly graceful graphs,}
\newblock {Discrete Math.,}
\newblock {29, 39-46, 1980.}
 https:// core.ac.uk/download/pdf/82643761.pdf
 
 \bibitem[Moulton, (1989)]{Mo}
 D. Moulton, 
 \textit{Graceful Labellings of Triangular Snakes,}
 \newblock {Ars. Combin.,}
 \newblock {28 }
 \newblock {2-13, 1989.}

 \bibitem[O'keefe, (1961)]{OK}
E.S. O'Keefe,
\textit{Verification of a conjecture of Th. Skolem,} 
\newblock {Math. Scand.,}
\newblock {9, 80-82, 1961.}

\bibitem[Rao & Sahoo, (2015)]{RS}
S. B. Rao and U. K. Sahoo,
\textit{Embeddings in Eulerian graceful graphs,} 
\newblock {Australasian J. Comb.,}
\newblock {62(1), 128-139, 2015.}

\bibitem[Rosa, (1967)]{R}
A. Rosa,
\textit{On certain valuations of the vertices of a graph,}
\newblock {Theory of Graphs (International Symposium , Rome, July 1966), Gordon and Breach, N. Y. and Dunod Paris,}
\newblock{349-355},
\newblock{1967.}

\bibitem[Sekar, (2002)]{S}
C. Sekar, 
\textit{Studies in Graph Theory,}
\newblock {Ph.D. Thesis, Madurai Kamaraj University,}
 \newblock{2002.}

\bibitem[Skolem, (1957)]{ST}
Th. Skolem,
\newblock {On certain distributions of integers in pairs with given differences,} 
\newblock {Math. Scand.,}
\newblock {5, suppl., 57-58, 1957.}

\bibitem[West, (2002)]{W}
D. B. West, 
\textit{Introduction to graph theory second edition,}
Pearson Education (Singapore) Pte. Ltd., Indian Braneh, 482 F.I.E. Patparganj,
Delhi 1 10 092, India, 2002.

\bibitem[Yang & Wang, (1992)]{YW}
Y. C. Yang and X. G. Wang, 
\textit{On the gracefulness of the product $C_{n}xP_{2}$,}
\newblock { J. Math. Research and Exposition, }
\newblock {1, 143-148, 1992.}

\end{thebibliography}
\end{document}